\documentclass[12pt, reqno]{amsart}

\usepackage[margin=1.2in]{geometry}
\usepackage{amsmath,amssymb,amsthm,mathtools}
\usepackage{booktabs,longtable,array,pdflscape}
\usepackage{needspace}
\usepackage[noadjust]{cite}
\usepackage{hyperref}
\hypersetup{
    colorlinks=true,
    citecolor=black,
    linkcolor=black,
    urlcolor=black
}

\theoremstyle{plain}
\newtheorem{theorem}{Theorem}[section]
\newtheorem{proposition}[theorem]{Proposition}
\newtheorem{lemma}[theorem]{Lemma}
\newtheorem{corollary}[theorem]{Corollary}

\theoremstyle{definition}
\newtheorem{definition}[theorem]{Definition}
\newtheorem{problem}[theorem]{Problem}

\theoremstyle{remark}
\newtheorem{remark}[theorem]{Remark}

\title[Spectral Equivalence of Complex Solvable Lie Algebras]
{Spectral Invariants of Complex Solvable Lie Algebras: Nilradical Weights and Hyperplane Arrangements}

\author{Gary Hu}
\email{gh7@williams.edu}

\author{Rongwei Yang}
\email{ryang@albany.edu}
\date{}


\begin{document}

\begin{abstract}
    A central open problem in Lie theory is the classification of finite-dimensional complex solvable Lie algebras, but classification up to isomorphism becomes increasingly difficult as the dimension grows and the number of non-isomorphic families explodes. This motivates a classification method by studying coarser spectral invariants arising from the characteristic polynomial of the adjoint representation. We prove that the characteristic polynomial is determined by the generalized weights of the adjoint action on the nilradical. This addresses the problem of expressing the spectral index in Lie-theoretic terms, gives sharp bounds in terms of the nilradical and quotient dimensions, and characterizes spectral equivalence. We then resolve a second problem concerning the higher Betti numbers of the eigenvariety complement. We endow the distinct non-$z_0$ factors of the characteristic polynomial with a natural matroid structure, which we call the spectral matroid, and use the Orlik--Solomon algebra of the associated hyperplane arrangement to determine the Betti numbers and Poincar\'{e} polynomial combinatorially and to prove log-concavity of the Betti sequence. Finally, we apply our theory to solvable Lie algebras with abelian nilradical to obtain explicit characteristic polynomials and spectral-equivalence criteria, including examples of non-isomorphic algebras with identical characteristic polynomials.
    \end{abstract}

\maketitle

\section{Introduction}

\subsection{Background}

This paper studies finite-dimensional Lie algebras over the complex field
$\mathbb{C}$. Classifying these algebras is one of the oldest
and most important problems in Lie theory. While the classification of
semisimple Lie algebras has been complete for over a century through
the works of Killing and Cartan \cite{Killing1888, Cartan1894}, the solvable
case remains open even in low dimensions such as eight or nine. The challenge is
largely due to the presence of a great many non-isomorphic classes of solvable
Lie algebras. In \cite{KeyYang2020}, Azari-Key and the second author introduced
the notion of spectral equivalence for Lie algebras. It is known that for finite-dimensional
complex simple Lie algebras, the notion of spectral equivalence agrees with isomorphism
\cite{Feng2026,GengLiuWang2024,Yang2024}. However, for solvable Lie algebras, spectral
equivalence is weaker (and therefore more flexible) than isomorphism. This leads to the
natural question whether classifying solvable Lie algebras under spectral equivalence
might be more fruitful. This paper develops the general theory needed to pursue that question.

\subsection{Solvable Lie Algebras and Nilradicals}

We first recall the decomposition theorems of Levi \cite{Levi1905} and
Malcev \cite{Malcev1945}. Levi's Theorem states that every finite-dimensional
Lie algebra $\mathfrak g$ decomposes as a semidirect product
$\mathfrak g=\mathfrak l\ltimes\operatorname{Rad}(\mathfrak g)$, where
$\operatorname{Rad}(\mathfrak g)$ is the unique maximal solvable ideal and
$\mathfrak l$ is a semisimple subalgebra. Malcev adds that the Levi factor
$\mathfrak l$ is unique up to inner automorphism. Together, these results
reduce the classification of finite-dimensional Lie algebras to the study of
solvable Lie algebras and the actions of semisimple Lie algebras on them. Of
these tasks, the classification of solvable Lie algebras $\mathfrak{g}$ is substantially more
challenging. The principal difficulty is therefore the classification of solvable
Lie algebras and of the compatible actions occurring in Levi extensions.

Since the radical $\operatorname{Rad}(\mathfrak g)$ is itself solvable, its
structure is closely tied to its maximal nilpotent ideal, the nilradical. Every
finite-dimensional solvable Lie algebra $\mathfrak g$ over $\mathbb{C}$ contains
a unique maximal nilpotent ideal $\operatorname{Nil}(\mathfrak g)$, and one has
$[\mathfrak g,\mathfrak g]\subseteq\operatorname{Nil}(\mathfrak g)$. After
choosing a vector-space complement $\mathfrak v$ to the nilradical, we may write
$\mathfrak g=\operatorname{Nil}(\mathfrak g)\oplus\mathfrak v$. This decomposition
is not canonical, and $\mathfrak v$ need not be a Lie subalgebra. The adjoint
action of elements of $\mathfrak v$ on $\operatorname{Nil}(\mathfrak g)$ gives
derivations of the nilradical, well-defined modulo inner derivations after
passing to the quotient $\mathfrak g/\operatorname{Nil}(\mathfrak g)$.

This idea motivates the classification strategy initiated by Mubarakzyanov
\cite{Mubarakzyanov1963_v5_1}, which first fixes a nilpotent Lie algebra and
then determines all possible solvable extensions having that algebra as
nilradical.

Classification efforts have also proceeded by dimension. Complete lists of
nilpotent Lie algebras are available up to dimension seven
\cite{Morozov1958,Seeley1993,Gong1998}, while partial and computer-assisted
classifications are available for several subclasses in dimensions eight and nine
\cite{Tsagas2000}. In the
solvable case, the classification is complete up to dimension seven. The
classification up to dimension six was completed by
\cite{Mubarakzyanov1963_v5_1,Mubarakzyanov1963_v5_2,Mubarakzyanov1963_v6,
Turkowski1990}, with \cite{SnoblWinternitz2014} providing an exposition. The
classification in dimension seven was finished by assembling partial results
organized according to the dimension of the nilradical: the remaining cases
were those with codimension-one nilradicals \cite{Parry2007}, four-dimensional
nilradicals \cite{HindelehThompson2008}, and five-dimensional nilradicals
\cite{Le2022}. The difficulty of extending these lists to higher dimensions
is primarily due to the rapid growth of non-isomorphic classes.

Following Mubarakzyanov's method, substantial progress has been made for
specific nilradicals, including abelian \cite{NdogmoWinternitz1994}, triangular
\cite{TremblayWinternitz1998}, and Borel nilradicals
\cite{WinternitzSnobl2012}, as well as filiform and quasifiliform nilradicals
\cite{SnoblWinternitz2005,SnoblWinternitz2009,Ancochea2006,
WangLinDeng2008,SnoblKarasek2010,Snobl2011,Ancochea2011}. Some nilpotent Lie
algebras, called characteristically nilpotent, cannot occur as nilradicals of
non-nilpotent solvable Lie algebras \cite{Khakimdzhanov1989,GozeHakimjanov1994}.

\needspace{8\baselineskip}
\subsection{Solvable Lie Algebras with Abelian Nilradical}

The fixed-nilradical strategy described above leads, in this case, to a
classification problem for commuting derivations and extension data.  Let
$\mathfrak{s}$ be a finite-dimensional complex solvable Lie algebra with
abelian nilradical $\mathfrak{n}$.  Then
$[\mathfrak{s},\mathfrak{s}]\subseteq\mathfrak{n}$ and
$\mathfrak{s}/\mathfrak{n}$ is abelian.  If
$\mathfrak{s}=\mathfrak{n}\oplus\mathfrak{v}$ as a vector space, the
restrictions $\operatorname{ad}(x)|_{\mathfrak{n}}$, $x\in\mathfrak{v}$,
commute.  Over $\mathbb{C}$ they are simultaneously upper-triangular, with
linear weights on the diagonal.  The remaining upper-triangular entries and
the brackets in $\mathfrak{v}$ contain the extension data.  Ndogmo and
Winternitz developed the structure theory, decomposability criteria, center
bounds, and examples for this class \cite{NdogmoWinternitz1994}; a broader
account of fixed-nilradical classification was given in
\cite{SnoblWinternitz2014}.

The class also has an invariant-theoretic role.  Patera, Sharp, Winternitz,
and Zassenhaus calculated generalized Casimir invariants for real Lie algebras
and discussed their role in representation theory and physics
\cite{PateraSharpWinternitzZassenhaus1976}.  Ndogmo and Winternitz constructed
generalized Casimir invariants for solvable Lie algebras with abelian
nilradicals \cite{NdogmoWinternitzCasimir1994}, and Ndogmo determined the
number of coadjoint invariant functions and the variables on which they can
depend \cite{Ndogmo2000}.  These works placed coadjoint orbits and their
invariants alongside the classification of the underlying Lie algebras.

Solvable Lie algebras also occur as symmetry algebras of differential
equations.  Their subalgebras provide invariant variables and group-invariant
solutions, making the identification of the abstract symmetry algebra a
basic step in Lie symmetry reduction \cite{Olver1993,NdogmoWinternitz1994}.
Ndogmo and
Winternitz emphasized the abelian case for constant-coefficient equations and
for equations reducible to that form \cite{NdogmoWinternitz1994}.  Solvable
symmetry algebras also appear in the construction of solutions for partially
integrable evolution equations \cite{Winternitz1990}.  Related fixed-nilradical
work connected generalized
Casimir invariants with the Petrov classification of Einstein spaces and with
the classification of Riemannian manifolds \cite{SnoblWinternitz2005}; that
example used a broader family of nilradicals.

At the group level, a simply connected Lie group whose Lie algebra has
abelian nilradical has a vector group as its nilpotent normal subgroup and an
abelian quotient, with the quotient action represented by commuting matrices.
This gives an
explicit setting for the representation theory and harmonic analysis of
solvable groups \cite{BakloutiFujiwaraLudwig2021}.  The corresponding real
groups also occur in solvmanifold and homogeneous geometry: Bock studied
lattices, formality, symplectic and K\"ahler structures, and the Hard Lefschetz
property for low-dimensional solvmanifolds \cite{Bock2016}, while Heber and
Lauret studied Einstein solvmanifolds and Ricci soliton solvmanifolds
\cite{Heber1998,Lauret2011}.  Within the abelian-nilradical class, contact
and Frobenius structures have been characterized in a semidirect subfamily
by Alvarez, Rodr\'iguez-Vallarte, and Salgado
\cite{AlvarezRodriguezVallarteSalgado2018}.

These examples also lead back to the algebraic questions studied in this
paper.  Ayupov and Khudoyberdiyeva treated local derivations in the
one-dimensional-complement and maximal-complement cases
\cite{AyupovKhudoyberdiyeva2021}.  In the present paper, the weights
determine the characteristic polynomial and the associated hyperplane
arrangement, while the extension data may distinguish non-isomorphic
algebras.

\subsection{Goals, Main Results, and Organization}

This paper is the first part of a series building on the spectral approach introduced in \cite{HuZhang2019,KeyYang2020}. Its first goal is general: for
every finite-dimensional complex solvable Lie algebra, we express the
characteristic polynomial in terms of the generalized weights of the action on
the nilradical. Its second goal is an application: for solvable algebras with
abelian nilradical, we translate those weights into the commuting-matrix
blocks of the Ndogmo--Winternitz normal form. The general results appear in
Sections~\ref{sec:general-theory}--\ref{sec:topological-invariants} and the
abelian specialization begins in Section~\ref{sec:abelian-nilradical}.
Subsequent papers apply the same framework to other structured nilradicals.

Azari-Key and the second author introduced spectral equivalence and the spectral index, and
obtained initial estimates and examples for these invariants
\cite{KeyYang2020}. The present paper develops the Lie-theoretic structure
behind them. Theorem~\ref{thm:nilradical-factorization} expresses the full
characteristic polynomial in terms of the nilradical weight multiset;
Theorem~\ref{thm:k-bounds} gives sharp dimension bounds for the spectral
index; and Proposition~\ref{prop:weight-equivalence} gives a criterion
for spectral equivalence in terms of the multiplicity of the factor $z_0$ and the
nonzero weights. The topological part is also new at the level of Lie
algebras: the higher cohomology of the eigenvariety complement is organized by
the matroid of the distinct factors other than $z_0$, leading to the
Poincar\'{e} polynomial, degree, bound, and log-concavity results in
Corollary~\ref{cor:poincare-matroid}, Theorem~\ref{thm:poincare-degree}, and
Corollaries~\ref{cor:betti-bound} and~\ref{cor:betti-log-concave}.
The exposition is organized as follows. Sections~\ref{sec:general-theory} and
\ref{sec:topological-invariants} contain the general theory. We define the
nilradical weight multiset, prove the characteristic-polynomial factorization,
derive the bounds for the spectral index, characterize spectral equivalence,
and identify the cohomology of the eigenvariety complement with the
Orlik--Solomon algebra of the associated spectral matroid. This gives the
degree, upper bound, and log-concavity results for the Betti numbers.

Section~\ref{sec:abelian-nilradical} contains the application. We derive the
block factorization in the Ndogmo--Winternitz normal form. We identify the
extension data that do not affect the characteristic polynomial, prove that
every finite spanning collection of nonzero weights occurs, and identify
the spectral matroid with the affine-dependence matroid of the distinct
nonzero weights. We then provide tables classifying the characteristic
polynomials and associated invariants for the normal form, the canonical
maximal-center families, and the eight-dimensional example in
\cite{NdogmoWinternitz1994}. Subsection~\ref{sec:abelian-worked-example}
provides an example showing how the matrix pencil, weight
multiset, spectral matroid, and extension data fit together.

\section{Characteristic Polynomials, Spectral Equivalence and Spectral Index}

The spectral theory of a matrix $A \in M_k(\mathbb{C})$ is traditionally developed via its characteristic polynomial $Q_A(\lambda) = \det(\lambda I - A)$. This invariant records the eigenvalues and algebraic multiplicities, and is one
of the basic invariants used to study the structure of a matrix. However, algebraic structures arising from physics and geometry are rarely generated by a single matrix, but rather by sets of non-commuting matrices $\{A_1, \dots, A_n\}$. Developing a spectral theory for such sets presents an immediate challenge, since non-commutativity prevents one from simultaneously diagonalizing the matrices or placing their spectral data into a single polynomial invariant using standard techniques.

\subsection{The Work of Dedekind and Frobenius}
The classical forerunner of a multivariate spectral theory is the theory of the group determinant, introduced by Dedekind and Frobenius in the 1890s. Although arising in a purely algebraic context, this construction displays the important features of a multivariate spectral invariant.

We first recall some basic notions in group representation theory. Let $G$ be a finite group of order $n+1$, written as $G = \{1, g_1, \ldots, g_n\}$. A unitary representation is a homomorphism $\pi: G\to U({\mathcal H})$, where ${\mathcal H}$ is a Hilbert space and $U({\mathcal H})$ denotes the group of unitary operators (or matrices) on ${\mathcal H}$. The dimension of the representation $\pi$, denoted by $d_\pi$, is the dimension of ${\mathcal H}$. Two unitary representations $\pi: G\to U({\mathcal H})$ and $\rho: G\to U({\mathcal K})$ are said to be unitarily equivalent, denoted by $\pi\cong\rho$, if there exists a unitary operator $V: {\mathcal H}\to {\mathcal K}$ such that $\pi (g)=V^*\rho(g)V$ for $g\in G$. Furthermore, a unitary representation is said to be irreducible if $\pi(G)$ has no nontrivial common reducing subspace in ${\mathcal H}$. We let $\widehat G$ denote the set of equivalence classes of irreducible unitary representations of $G$ (called the unitary dual of $G$).
Finally, the left regular representation $\lambda_G: G\to U(\ell^2(G))$ is defined by $\lambda_G(g)f(x)=f(g^{-1}x)$ for $x\in G$.

The group determinant with respect to a finite-dimensional representation $\pi$ is defined by
\[Q_\pi(z_0,\mathbf{z}) := \det\left(z_0 I + z_1 \pi(g_1) + \cdots + z_n \pi(g_n)\right),\]
where $\mathbf{z}=(z_1,\ldots,z_n)\in \mathbb{C}^n$.
Frobenius's breakthrough was the discovery that the factorization of this determinant for the left regular representation encodes the irreducible unitary representations of $G$.

\begin{theorem}[\cite{Dedekind1969,Frobenius1896,Dickson1902}]\label{thm:group-determinant}
    For any finite group $G=\{1,g_1,\dots,g_n\}$, we have
    \[
    Q_{\lambda_G}(z_0,\mathbf{z})=\prod_{[\pi]\in\widehat G} Q_\pi(z_0,\mathbf{z})^{d_\pi}.
    \]
    Moreover, every $Q_\pi$ is an irreducible polynomial, and $\pi\cong \rho$ if and only if $Q_\pi=Q_\rho$.
\end{theorem}

In fact, this theorem marks the beginning of representation theory for finite groups. For historical developments, see \cite{Curtis1992, Curtis1999, Dickson1921, Dickson1902, FormanekSibley1991}.

\subsection{Characteristic Polynomial}
Despite this early success, it took more than a century before the group determinant was extended to other algebraic settings. We briefly outline this history; a comprehensive survey is available in \cite{Yang2024}. In 2009, \cite{Yang2009} extended this notion to unital Banach algebras, initiating a line of research that has since attracted considerable attention \cite{CadeYang2013, Chagouel2015, DouglasYang2018, GoldbergYang2022, GrigorchukYang2017, HeWangYang2017, MaoQiaoWang2017, StessinYangZhu2011}. Motivated by these results, \cite{HuZhang2019} defined the characteristic polynomial for finite-dimensional Lie algebras.

\begin{definition} \label{def:char-poly}
Let $L$ be a finite-dimensional Lie algebra with basis $\{x_1, \dots, x_n\}$. The \emph{characteristic polynomial} of $L$ is
\[
Q_L(z_0, \mathbf{z}) := \det \bigg( z_0 I + \sum_{i=1}^n z_i \operatorname{ad} x_i \bigg),\qquad \mathbf{z}=(z_1,\ldots,z_n)\in \mathbb{C}^{n}.\]
\end{definition}

Here $\operatorname{ad}(x)$ denotes the adjoint map $y\mapsto[x,y]$.
The matrix inside this determinant is a matrix pencil: its entries depend
linearly on the variables $z_0,z_1,\ldots,z_n$.

This notion has already been applied to study simple \cite{GengLiuWang2024}, classical \cite{Feng2026}, and solvable Lie algebras \cite{HuZhang2019, Muller}. The solvable case is particularly interesting: by Lie's theorem, the adjoint operators can be simultaneously upper triangularized over $\mathbb{C}$, and therefore the characteristic polynomial factors completely into linear forms. Furthermore, the converse is also true.

\begin{theorem}[\cite{HuZhang2019,Yang2024}]
\label{thm:solvable-factorization-old}
A finite-dimensional complex Lie algebra $L$ is solvable if and only if its characteristic
polynomial $Q_L(z_0,\mathbf{z})$ factors as
\[
Q_L(z_0,\mathbf{z})=\prod_{j=1}^{\dim L}(z_0+\ell_j(\mathbf{z})),
\]
where each $\ell_j$ is a linear form in $\mathbf{z}=(z_1,\ldots,z_n)$.
\end{theorem}

\begin{theorem}[\cite{KeyYang2020}]\label{thm:nilpotent-char-poly}
    A finite-dimensional complex Lie algebra $L$ is nilpotent if and only if
    \[
    Q_L(z_0, \mathbf{z})=z_0^{\dim L}.
    \]
\end{theorem}

\subsection{Spectral Equivalence and Spectral Index}
One motivation for developing a theory of characteristic polynomials for Lie algebras is the expectation that these polynomials will give spectral invariants capable of distinguishing and classifying Lie algebras. Observe that if two Lie algebras $L_1$ and $L_2$ are isomorphic, any isomorphism corresponds to a change of basis matrix $B \in \operatorname{GL}_n(\mathbb{C})$. A direct calculation shows that their characteristic polynomials must satisfy $Q_{L_2}(z_0, \mathbf{z}) = Q_{L_1}(z_0, \mathbf{z} B)$. Motivated by this observation, \cite{KeyYang2020} introduced a coarser classification criterion:

\begin{definition}
    Two $n$-dimensional complex Lie algebras $L_1$ and $L_2$ are
    \emph{spectrally equivalent} over $\mathbb{C}$, denoted $L_1\sim L_2$, if
    \[
    Q_{L_2}(z_0,\mathbf{z})=Q_{L_1}(z_0,\mathbf{z} B)
    \]
    for some $B\in\operatorname{GL}_n(\mathbb{C})$.
\end{definition}

We will primarily study this relation for solvable Lie algebras.

Spectral equivalence is strictly weaker than isomorphism: if $L_1 \cong L_2$, then $L_1 \sim L_2$, but the converse fails in general. However, the converse does hold in certain important cases, such as for complex simple Lie algebras:

\begin{theorem}[\cite{Yang2024}, Corollary 3.29]
    Let $L_1$ and $L_2$ be two finite-dimensional complex simple Lie algebras. Then they are isomorphic if and only if they are spectrally equivalent.
\end{theorem}

This does not extend to the broader class of solvable Lie algebras easily. Nevertheless, the factorization from Theorem \ref{thm:solvable-factorization-old} still gives rise to an interesting invariant, introduced in \cite{KeyYang2020}, which can distinguish some non-isomorphic solvable Lie algebras.

\begin{definition}
    Let $L$ be a finite-dimensional complex Lie algebra. Its \emph{spectral index} $\kappa(L)$ is defined as the number of distinct irreducible factors in the factorization of $Q_L(z_0, \mathbf{z})$.
\end{definition}

It is not hard to see that the spectral index is invariant under spectral
equivalence. For solvable Lie algebras, it counts the number of distinct
linear factors of $Q_L$. This index has been defined for arbitrary tuples of
matrices and investigated extensively in \cite{StessinYang2026}. However, for
Lie algebras, it remains quite mysterious, despite having been computed in
numerous examples and partial results. For a square matrix $A$, write
$\sigma(A)$ for its spectrum. Then $|\sigma(A)|$ is the number of distinct
eigenvalues of $A$. The following lower bound for $\kappa(L)$ is known.
\begin{proposition}[\cite{KeyYang2020}]
Let $L$ be a solvable Lie algebra with basis
$\{x_1,x_2,\dots,x_n\}$. Then $\kappa(L)\geq \max\{ |\sigma(\operatorname{ad} x_i)|: 1\leq i\leq n\}$.
\end{proposition}

The paper \cite{KeyYang2020} also raised the following question, aiming to better characterize $\kappa(L)$. 

\begin{problem}\label{prob:k-invariant}
    Can $\kappa(L)$ be expressed in terms of other Lie-theoretic invariants?
\end{problem}

We return to this question after introducing the nilradical weight multiset
and proving the general factorization theorem.

\section{Spectral Index and Nilradical Weight Multiset} \label{sec:general-theory}

The results in this section are valid for every finite-dimensional complex
solvable Lie algebra. We now identify the weights of the action on the
nilradical and use them to factor the characteristic polynomial.

\subsection{The Spectral Index \texorpdfstring{$\kappa(L)$}{kappa(L)}}

Let $\mathfrak{s}$ be a finite-dimensional solvable Lie algebra over $\mathbb{C}$ with nilradical $\mathfrak{n}$. Choose a vector-space complement $\mathfrak{v}$ so that $\mathfrak{s}=\mathfrak{n}\oplus \mathfrak{v}$, and write $k=\dim \mathfrak{v}$. Since $\mathfrak{n}$ is an ideal, the adjoint action restricts to a representation $\rho:\mathfrak{s}\to\mathfrak{gl}(\mathfrak{n})$, where $\rho(y)=\operatorname{ad}(y)|_{\mathfrak{n}}$.

By Lie's theorem, there is a basis of $\mathfrak{n}$ in which every $\rho(y)$ is upper triangular. Thus there are linear functionals $\lambda_1,\dots,\lambda_M\in\mathfrak{s}^*$, where $M=\dim\mathfrak{n}$, such that the diagonal entries of $\rho(y)$ are $\lambda_1(y),\dots,\lambda_M(y)$ for every $y\in\mathfrak{s}$. If $x\in\mathfrak{n}$, then $\rho(x)$ is nilpotent, so every $\lambda_i$ vanishes on $\mathfrak{n}$. Hence each $\lambda_i$ descends to $(\mathfrak{s}/\mathfrak{n})^*$, and after choosing $\mathfrak{v}$, we regard it as an element of $\mathfrak{v}^*$.

\begin{definition}\label{def:nilradical-weights}
The \emph{nilradical weight multiset} of $\mathfrak{s}$ is $\Delta_{\mathfrak{s}}:=\{\lambda_1|_{\mathfrak{v}},\dots,\lambda_M|_{\mathfrak{v}}\}$, counting multiplicity. We write $\Delta$ for the corresponding set of distinct weights. If $\alpha\in\Delta$, let $m_\alpha$ denote the multiplicity of $\alpha$.
\end{definition}

\begin{lemma}\label{lem:weights-well-defined}
The nilradical weight multiset $\Delta_{\mathfrak s}$ is independent of the
choice of triangularizing basis.
\end{lemma}

\begin{proof}
Choose a composition series
$0=N_0\subset N_1\subset\cdots\subset N_M=\mathfrak n$ of the
$\mathfrak{s}$-module $\mathfrak n$. Each one-dimensional quotient
$N_i/N_{i-1}$ is acted on by a linear functional $\lambda_i\in\mathfrak{s}^*$.
By the Jordan--Hölder theorem, the multiset of these composition factors, and
hence the multiset of the $\lambda_i$, is independent of the chosen series.
For $x\in\mathfrak n$, the operator $\operatorname{ad}(x)|_{\mathfrak n}$ is
nilpotent, so its induced scalar on every one-dimensional quotient is zero.
Thus each $\lambda_i$ vanishes on $\mathfrak n$ and descends to
$(\mathfrak{s}/\mathfrak n)^*$. Choosing $\mathfrak v$ identifies this dual
space with $\mathfrak v^*$, which proves the claim.
\end{proof}

\begin{remark}
The weights in Definition~\ref{def:nilradical-weights} are the diagonal
functionals from the triangular form. No simultaneous diagonalization is
assumed, and the operators may have nonzero strictly upper-triangular parts.
\end{remark}

Choose a basis $\mathcal{B}=\{x_1,\dots,x_M,f_1,\dots,f_k\}$ for $\mathfrak{s}$ adapted to the decomposition $\mathfrak{s}=\mathfrak{n}\oplus \mathfrak{v}$, where $M = \dim \mathfrak{n}$ and $k = \dim \mathfrak{v}$. Let $N = M + k$ denote the total dimension of $\mathfrak{s}$. We associate this basis with the set of indeterminates
$\mathbf{z}=(\mathbf{z}_{\mathfrak n},\mathbf{z}_{\mathfrak v}):=(z_{x_1},\dots,z_{x_M},z_{f_1},\dots,z_{f_k})$. For each weight $\alpha\in\Delta$, define the linear form
\[\ell_\alpha(\mathbf{z}_{\mathfrak v}):=\sum_{j=1}^k \alpha(f_j)\,z_{f_j}.\]
In these adapted coordinates, the characteristic polynomial is
$Q_{\mathfrak{s}}(z_0,\mathbf{z})=\det \mathbf A(z_0,\mathbf{z})$, where
\[
\mathbf{A}(z_0, \mathbf{z}):=z_0 I_N + \sum_{i=1}^M z_{x_i}\operatorname{ad}(x_i) + \sum_{j=1}^k z_{f_j}\operatorname{ad}(f_j),
\]
and $I_N$ is the identity matrix of size $N\times N$.

\subsection{A Factorization of the Characteristic Polynomial}

With the weight structure established, we use it to study the characteristic polynomial. Our approach is guided by the goal of obtaining a statement analogous to the Dedekind--Frobenius Theorem (Theorem \ref{thm:group-determinant}), where the factorization of the group determinant is indexed by the unitary dual of a finite group. In the Lie algebra setting, it is natural to ask whether the characteristic polynomial $Q_\mathfrak{s}(z_0,\mathbf{z})$ admits a similar structural interpretation.

For finite-dimensional simple Lie algebras, this question has been addressed through the lens of invariant theory and Weyl group actions. Feng, Liu, and Wang \cite{Feng2026} showed that for classical Lie algebras and the exceptional algebra $G_2$, the characteristic polynomial associated with a finite-dimensional representation decomposes into a product of irreducible orbital factors. These factors are uniquely determined by the orbits of the weights of the representation under the action of the corresponding Weyl group.

For solvable Lie algebras, the generalized weights of the action on the nilradical determine the linear factors of the characteristic polynomial. Let $\Delta\subset \mathfrak{v}^*$ be the set of distinct weights from Definition~\ref{def:nilradical-weights}, and recall that $m_\alpha$ denotes the multiplicity of $\alpha$ in the nilradical weight multiset.

\begin{theorem}\label{thm:nilradical-factorization}
Let $\mathfrak{s}=\mathfrak{n}\oplus \mathfrak{v}$ be solvable. Then
\[
Q_{\mathfrak{s}}(z_0,\mathbf{z})=z_0^{\dim \mathfrak{v}}\prod_{\alpha\in\Delta}\bigl(z_0+\ell_\alpha(\mathbf{z}_{\mathfrak v})\bigr)^{m_\alpha}.
\]
Consequently, the spectral index is $\kappa(\mathfrak{s})=|\Delta\cup\{0\}|$.
\end{theorem}

\begin{proof}

With respect to the decomposition $\mathfrak{s}=\mathfrak{n}\oplus \mathfrak{v}$, the matrix pencil $\mathbf A(z_0,\mathbf{z})$ is block upper-triangular:
\[
\mathbf A(z_0,\mathbf{z})=
\begin{pmatrix}
\mathbf A_{\mathfrak n}(z_0,\mathbf{z}) & *\\
0 & \mathbf A_{\mathfrak{s}/\mathfrak n}(z_0,\mathbf{z})
\end{pmatrix}.
\]
Here $\mathbf A_{\mathfrak n}$ is the restriction of the pencil to the invariant subspace $\mathfrak n$, while $\mathbf A_{\mathfrak{s}/\mathfrak n}$ denotes the induced pencil on the quotient $\mathfrak{s}/\mathfrak n$, identified with $\mathfrak{v}$ as a vector space. Taking determinants gives
\[
Q_{\mathfrak{s}}(z_0,\mathbf{z})
=
\det(\mathbf A_{\mathfrak n}(z_0,\mathbf{z}))\,
\det(\mathbf A_{\mathfrak{s}/\mathfrak n}(z_0,\mathbf{z})).
\]

We first compute the quotient block. Since $\mathfrak{s}$ is solvable, Lie's theorem gives a basis in which all operators $\operatorname{ad}(\mathfrak{s})$ are upper triangular. Hence all operators $\operatorname{ad}([\mathfrak{s},\mathfrak{s}])$ are strictly upper triangular, and therefore nilpotent. By Engel's theorem, $[\mathfrak{s},\mathfrak{s}]$ is a nilpotent ideal, so by maximality of the nilradical,
\[
[\mathfrak{s},\mathfrak{s}]\subseteq \mathfrak n.
\]
Thus the induced adjoint action of $\mathfrak{s}$ on $\mathfrak{s}/\mathfrak n$ is trivial. Hence
\[
\mathbf A_{\mathfrak{s}/\mathfrak n}(z_0,\mathbf{z})=z_0I_{\dim \mathfrak{v}},
\qquad
\det(\mathbf A_{\mathfrak{s}/\mathfrak n}(z_0,\mathbf{z}))=z_0^{\dim \mathfrak{v}}.
\]

Next, consider $\det(\mathbf A(z_0,\mathbf{z})|_{\mathfrak n})$. Choose a basis
of $\mathfrak n$ in which all operators
$\rho(y)=\operatorname{ad}(y)|_{\mathfrak n}$, $y\in\mathfrak s$, are upper
triangular. Let $\lambda_1,\dots,\lambda_M$ be the corresponding weights. Since
each $\lambda_i$ vanishes on $\mathfrak n$, the variables
$z_{x_1},\dots,z_{x_M}$ do not appear on the diagonal of the nilradical block.
If $\lambda_i|_{\mathfrak{v}}=\alpha$, then the corresponding diagonal entry in the matrix $\mathbf A_{\mathfrak{n}}(z_0,\mathbf{z})$ is
$z_0+\ell_\alpha(\mathbf{z}_{\mathfrak v})$. Hence the determinant of this nilradical block contributes $\prod_{\alpha\in\Delta}(z_0+\ell_\alpha(\mathbf{z}_{\mathfrak v}))^{m_\alpha}$.

Combining this with the quotient block gives
\[
Q_{\mathfrak{s}}(z_0,\mathbf{z})=z_0^{\dim \mathfrak{v}}\prod_{\alpha\in\Delta}(z_0+\ell_\alpha(\mathbf{z}_{\mathfrak v}))^{m_\alpha}.
\]
The distinct linear factors are $z_0$ together with the factors $z_0+\ell_\alpha(\mathbf{z}_{\mathfrak v})$ for those $\alpha\in\Delta$ with $\alpha\neq 0$. Equivalently, the distinct factors are indexed by $\Delta\cup\{0\}$, where the element $0$ corresponds to the factor $z_0$. Hence \[ \kappa(\mathfrak{s})=|\Delta\cup\{0\}|. \]
\end{proof}

This result resolves Problem \ref{prob:k-invariant}.

\begin{remark}
    For solvable Lie algebras, Theorem~\ref{thm:nilradical-factorization} recovers the nilpotence criterion in Theorem~\ref{thm:nilpotent-char-poly}. If $L$ is nilpotent, then $L=\mathfrak{n}$ and $Q_L(z_0,\mathbf{z})=z_0^{\dim L}$. Conversely, if $Q_L(z_0,\mathbf{z})=z_0^{\dim L}$, then every nilradical weight is zero. Hence every $f\in \mathfrak{v}$ acts nilpotently on $\mathfrak{n}$, and the induced action on $L/\mathfrak{n}$ is trivial. Thus $\operatorname{ad}(f)$ is nilpotent on $L$. Since, for a finite-dimensional solvable Lie algebra over $\mathbb{C}$, the nilradical is the set of ad-nilpotent elements \cite[Ch. 1, Sect. VI, Theorem 11]{GozeHakimjanov1996}, we get $f\in\mathfrak{n}$. Therefore $\mathfrak{v}=0$ and $L=\mathfrak{n}$ is nilpotent.
\end{remark}

We can use this theorem to obtain sharp bounds for the spectral index $\kappa(\mathfrak{s})$.

\begin{theorem}\label{thm:k-bounds}
Consider a solvable Lie algebra $\mathfrak{s}=\mathfrak{n}\oplus \mathfrak{v}$. Then we have
\[
\dim \mathfrak{v}+1\leq \kappa(\mathfrak{s})\leq \dim\mathfrak{n}+1.
\]
Moreover, $\kappa(\mathfrak{s})=\dim \mathfrak{v}+1$ if and only if the nonzero weights $\Delta\setminus\{0\}$ form a basis of $\mathfrak{v}^*$. Also, $\kappa(\mathfrak{s})=\dim\mathfrak{n}+1$ if and only if $0\notin\Delta$ and every distinct weight has multiplicity one, that is, $m_\alpha=1$ for all $\alpha\in\Delta$.
\end{theorem}

\begin{proof}
    We prove the two bounds separately.

    \emph{Lower bound.} We first show that the set of weights spans the dual space. Suppose $\operatorname{span}(\Delta)\subsetneq \mathfrak{v}^*$. Then there exists a nonzero $f_0\in \mathfrak{v}$ such that $\alpha(f_0)=0$ for all $\alpha\in\Delta$. The eigenvalues of $\operatorname{ad}(f_0)|_{\mathfrak{n}}$ are $\{\alpha(f_0)\mid\alpha\in\Delta\}=\{0\}$, so this operator is nilpotent. Moreover, as shown in the proof of Theorem~\ref{thm:nilradical-factorization}, the induced action on $\mathfrak{s}/\mathfrak{n}$ is trivial, hence nilpotent. Therefore $\operatorname{ad}(f_0)$ is globally nilpotent on $\mathfrak{s}$, so $f_0$ belongs to the nilradical $\mathfrak{n}$. This contradicts $f_0\in \mathfrak{v}\setminus\{0\}$.

    Thus $\operatorname{span}(\Delta)= \mathfrak{v}^*$. Since the zero vector does not contribute to the span, the set of nonzero weights $\Delta' := \Delta \setminus \{0\}$ must span $ \mathfrak{v}^*$. Consequently, $|\Delta'| \ge \dim( \mathfrak{v})$.

     Recall that $\kappa(\mathfrak{s})$ counts the number of distinct linear factors in $Q_{\mathfrak{s}}(z_0,\mathbf{z})$. By Theorem~\ref{thm:nilradical-factorization}, the extension $ \mathfrak{v}$ contributes the factor $z_0^{\dim  \mathfrak{v}}$. Thus the factor $z_0$ is always present in $Q_{\mathfrak{s}}$, coming either from the quotient block or, if $0\in\Delta$, also from zero nilradical weights. Hence the distinct factors are indexed by $\Delta\cup\{0\}$. We can compute the cardinality as:
        \[
        \kappa(\mathfrak{s}) = |\Delta \cup \{0\}| = |\Delta'| + 1.
        \]

    Combining this with the span inequality $|\Delta'| \ge \dim( \mathfrak{v})$, we obtain:
        \[
        \kappa(\mathfrak{s}) \ge \dim( \mathfrak{v}) + 1.
        \]

        Equality holds if and only if $|\Delta'| = \dim( \mathfrak{v})$, which occurs if and only if the nonzero weights $\Delta'$ form a basis of $ \mathfrak{v}^*$.

    \emph{Upper bound.} Since the nilradical weight multiset has total cardinality $\dim\mathfrak n$, we have $\sum_{\alpha\in\Delta}m_\alpha=\dim\mathfrak n$. In particular, $|\Delta|\leq\dim\mathfrak{n}$. If $0\in\Delta$, then $\kappa(\mathfrak{s})=|\Delta|\leq\dim\mathfrak{n}$. If $0\notin\Delta$, then $\kappa(\mathfrak{s})=|\Delta|+1\leq\dim\mathfrak{n}+1$. This proves the bound. Equality requires $0\notin\Delta$ and $|\Delta|=\dim\mathfrak{n}$. Since $\sum_{\alpha\in\Delta}m_\alpha=\dim\mathfrak{n}$ and every $m_\alpha\geq 1$, this happens if and only if $m_\alpha=1$ for every $\alpha\in\Delta$.
\end{proof}

\subsection{A Characterization of Spectral Equivalence}
The factorization of Theorem~\ref{thm:nilradical-factorization} allows the notion of spectral equivalence to be reformulated as a condition on the collection of weights. We denote by $\Delta_{\mathfrak s}^{\times}$ the multiset of nonzero nilradical weights of $\mathfrak{s}$ and by $e_0(\mathfrak{s})$ the total multiplicity of the factor $z_0$ in $Q_{\mathfrak{s}}$. Two multisets of vectors $E$ and $E'$ are said to be linearly equivalent, denoted $E\cong E'$, if there is a linear isomorphism $T:\operatorname{span} E\to\operatorname{span} E'$ with $T(E)=E'$, where multiplicities are preserved.

\begin{proposition}\label{prop:weight-equivalence}
Two finite-dimensional solvable Lie algebras $\mathfrak{s}$ and $\mathfrak{s}'$
are spectrally equivalent if and only if
$e_0(\mathfrak{s})=e_0(\mathfrak{s}')$ and the multisets
$\Delta_{\mathfrak s}^{\times}$ and $\Delta_{\mathfrak s'}^{\times}$ are
linearly equivalent, with multiplicities preserved.
\end{proposition}

\begin{proof}
We need to show that solvable Lie algebras $\mathfrak{s}$ and $\mathfrak{s}'$ are spectrally equivalent if and only if $e_0(\mathfrak{s})=e_0(\mathfrak{s}')$ and $\Delta_{\mathfrak s}^{\times}$ and $\Delta_{\mathfrak s'}^{\times}$ are linearly equivalent.

By Theorem~\ref{thm:nilradical-factorization}, the characteristic polynomial of
$\mathfrak s$ has the form
\[
Q_{\mathfrak s}
=
z_0^{e_0(\mathfrak s)}
\prod_{\alpha\in\Delta_{\mathfrak s}^{\times}}(z_0+\ell_\alpha),
\]
where the product is over the nonzero nilradical weight multiset, counting
multiplicity. Here each $\ell_\alpha$ is viewed as a linear form on the full
space of Lie-algebra variables by extending it by zero on the nilradical
variables. The same description holds for $\mathfrak s'$.

Suppose first that $\mathfrak s\sim\mathfrak s'$. Then
$Q_{\mathfrak s'}(z_0,\mathbf{z})=Q_{\mathfrak s}(z_0,\mathbf{z} B)$ for some
$B\in\operatorname{GL}_N(\mathbb{C})$. This change of variables fixes $z_0$.
Therefore the irreducible factor $z_0$ is carried to $z_0$, and unique
factorization gives $e_0(\mathfrak s)=e_0(\mathfrak s')$. The remaining factors
are the factors $z_0+\ell_\alpha$ associated with nonzero weights. Since an
invertible linear change of the Lie-algebra variables carries these factors to
factors other than $z_0$ and preserves multiplicity, it carries the factors of
$Q_{\mathfrak s}$ other than $z_0$ to the corresponding factors of
$Q_{\mathfrak s'}$. After subtracting the
common $z_0$ term, we obtain a linear isomorphism from
$\operatorname{span}\Delta_{\mathfrak s}^{\times}$ to
$\operatorname{span}\Delta_{\mathfrak s'}^{\times}$ carrying one multiset to
the other.

Conversely, suppose that $e_0(\mathfrak s)=e_0(\mathfrak s')$ and that such a
linear isomorphism
$T:\operatorname{span}\Delta_{\mathfrak s}^{\times}\to
\operatorname{span}\Delta_{\mathfrak s'}^{\times}$ exists. Let
$\mathfrak v$ and $\mathfrak v'$ be the complements used to represent the two
weight multisets. By the spanning argument in Theorem~\ref{thm:k-bounds},
$\operatorname{span}\Delta_{\mathfrak s}^{\times}=\mathfrak v^*$ and
$\operatorname{span}\Delta_{\mathfrak s'}^{\times}=(\mathfrak v')^*$, so
$T$ is already an isomorphism between the dual complement spaces. The equality
of $e_0$ and preservation of multiplicities imply that the two Lie algebras
have the same dimension. Choose complementary subspaces $U$ and $U'$ so that $\mathfrak s^*=\mathfrak v^*\oplus U$ and $(\mathfrak s')^*=(\mathfrak v')^*\oplus U'.$ It follows that $\dim U=\dim U'$. Choose any linear
isomorphism $S:U\to U'$ and define
$\widetilde T:=T\oplus S:\mathfrak s^*\to(\mathfrak s')^*$. Extending
$\widetilde T$ to the characteristic-polynomial variables by fixing $z_0$
carries each factor $z_0+\ell_\alpha$ to the corresponding factor for
$\mathfrak s'$, with the same multiplicity, and the exponent of $z_0$ also
agrees. Hence the characteristic polynomials are related by an invertible
change of variables, so $\mathfrak s\sim\mathfrak s'$.
\end{proof}

\begin{remark}
    If one works inside a family where $\dim \mathfrak{v}$ and the multiplicity of the zero
nilradical weight are fixed, this is equivalent to the simpler condition that
the full nilradical weight multisets are equivalent under
$\operatorname{GL}(\mathfrak{v}^*)$.
\end{remark}

Many classification lists of Lie algebras contain families of solvable ones depending
on continuous parameters. Let $\mathfrak{s}(\mathbf p)$ be such a family of
$N$-dimensional complex solvable Lie algebras, with parameter
$\mathbf p\in\mathcal P$. Since all nilpotent Lie algebras of a fixed dimension
have characteristic polynomial $z_0^N$, spectral equivalence is very coarse on
nilpotent algebras. For non-nilpotent solvable algebras, however, the nonzero
nilradical weights can vary with the parameters. This motivates the following definition.

\begin{definition}
    A parameterized family of Lie algebras $\mathfrak{s}(\mathbf{p})$ is \emph{spectrally rigid} if for any $\mathbf{p}, \mathbf{p}' \in \mathcal{P}$, the condition $\mathfrak{s}(\mathbf{p}) \sim \mathfrak{s}(\mathbf{p}')$ implies $\mathbf{p} = \mathbf{p}'$.
\end{definition}

By Proposition~\ref{prop:weight-equivalence}, spectral equivalence is controlled by the total multiplicity of the factor $z_0$ and by the orbit of the nonzero nilradical weight multiset. Thus, in any family where $e_0(\mathfrak{s}(\mathbf p))$ is constant and the nonzero weights span a fixed vector space $\mathfrak{v}^*$, spectral rigidity is equivalent to separation of the nonzero weight multisets under $\operatorname{GL}(\mathfrak{v}^*)$. For a parameterized family of Lie algebras $\mathfrak{s}(\mathbf{p})$, we let $\Delta^\times(\mathbf p)$ denote the nonzero nilradical weight multiset.
In fact, under the hypotheses below, Proposition~\ref{prop:weight-equivalence}
shows that two members are spectrally equivalent precisely when an element of
$\operatorname{GL}(\mathfrak v^*)$ carries one of these multisets to the other.
This gives the following criterion.

\begin{corollary}\label{cor:weight-rigidity}
Let $\mathfrak{s}(\mathbf p)$ be a family of solvable Lie algebras for which $e_0(\mathfrak{s}(\mathbf p))$ is constant and all nonzero weight multisets are viewed in a fixed space $\mathfrak{v}^*$. Then the family is spectrally rigid if and only if, for every $\mathbf p\neq\mathbf p'$, there is no $T\in\operatorname{GL}(\mathfrak{v}^*)$ such that $T(\Delta^\times(\mathbf p))=\Delta^\times(\mathbf p')$ as multisets.
\end{corollary}

\section{Topological Spectral Invariants} \label{sec:topological-invariants}

The main goal of this section is to use this observation
to attach topological invariants to hyperplane arrangements arising from the characteristic polynomial of $\mathfrak s$, and describe this theory in a combinatorial way. Theorem~\ref{thm:nilradical-factorization} expresses the characteristic polynomial $Q_{\mathfrak s}$ of a finite-dimensional complex solvable Lie algebra $\mathfrak s$ as a product of linear forms. Since $z_0$ is always a factor, we remove it and consider the distinct remaining factors $z_0+\ell_\alpha(\mathbf{z}_{\mathfrak v})$. These factors define a hyperplane arrangement in $\mathbb{C}^{\dim\mathfrak v+1}$, whose structure can be studied combinatorially through a matroid, which we call the spectral matroid $\mathcal M_{\mathfrak s}$. We show that the
cohomology ring of the eigenvariety complement $(V_{\mathfrak s})^c$ is
isomorphic to the Orlik--Solomon algebra of the spectral matroid
$\mathcal M_{\mathfrak s}$, so the Betti numbers of $\mathfrak s$ and its
Poincar\'{e} polynomial are determined by the Whitney numbers of the first
kind of $\mathcal M_{\mathfrak s}$, thereby providing an answer to
Problem~\ref{prob:higher-betti}.

\subsection{Eigenvarieties and Betti Numbers}
Let $\mathfrak s$ be a finite-dimensional complex solvable Lie algebra with
basis $\{x_1,\ldots,x_N\}$. At $z_0=0$, the matrix in
Definition~\ref{def:char-poly} is $\operatorname{ad}(\sum_i z_i x_i)$, which
has $\sum_i z_i x_i$ in its kernel. Hence \[Q_{\mathfrak{s}}(0,{\bf z})=\det (z_1\operatorname{ad}x_1+\cdots+z_n\operatorname{ad}x_n)=0, \ \ \ {\bf z}\in \mathbb{C}^N,\]
showing that $z_0$ always divides
$Q_{\mathfrak s}(z_0,\mathbf{z})$. Let $e_0(\mathfrak s)$ denote the total
multiplicity of the factor $z_0$ in $Q_{\mathfrak s}$, and set
$a=e_0(\mathfrak s)$. Write
\[
Q_{\mathfrak s}(z_0,\mathbf{z})
=:z_0^a p_{\mathfrak s}(z_0,\mathbf{z})
=z_0^a\prod_{\alpha\in\Delta\setminus\{0\}}
\bigl(z_0+\ell_\alpha(\mathbf{z}_{\mathfrak v})\bigr)^{m_\alpha},
\]
where $a\geq1$ and $z_0\nmid p_{\mathfrak s}$.

\begin{definition}\label{def:eigenvariety}
The \emph{eigenvariety} of $\mathfrak s$ is the zero locus of the polynomial $p_{\mathfrak s}$:
\[
V_{\mathfrak s}:=\{(z_0,\mathbf{z})\in\mathbb{C}^{\dim\mathfrak s+1}:
p_{\mathfrak s}(z_0,\mathbf{z})=0\},
\qquad
V_{\mathfrak s}^c:=\mathbb{C}^{\dim\mathfrak s+1}\setminus V_{\mathfrak s}.
\]
\end{definition}

\begin{definition}
Let $\mathfrak s$ be a finite-dimensional solvable Lie algebra. For each
$j\geq0$, the $j$-th \emph{Betti number} of $\mathfrak s$ is
$b_j(\mathfrak s):=\dim_{\mathbb{C}}H^j(V_{\mathfrak s}^c,\mathbb{C})$, where
$H^j$ denotes the $j$-th de Rham cohomology group with complex coefficients.
The \emph{Poincar\'{e} polynomial} of $\mathfrak s$ is
$P_{\mathfrak s}(t):=\sum_{j\geq0}b_j(\mathfrak s)t^j$.
\end{definition}

Some preliminary results are shown in \cite{KeyYang2020}.
\begin{proposition}\label{prop:first-betti}
Let $\mathfrak s$ be a finite-dimensional solvable Lie algebra. Then
\begin{itemize}
    \item $b_1(\mathfrak s)=\kappa(\mathfrak s)-1$.
    \item $\deg P_{\mathfrak s}\leq
    \dim\bigl(\mathfrak s/\operatorname{Nil}(\mathfrak s)\bigr)+1$.
\end{itemize}
\end{proposition}

However, the structure of the higher Betti numbers is posed as an open problem:
\begin{problem}\label{prob:higher-betti}
Can the Betti numbers $b_j(\mathfrak s)$, $j\geq2$, be described by known
invariants of $\mathfrak s$?
\end{problem}

The rest of this section answers this question. The cohomology is identified
with an Orlik--Solomon algebra, and the resulting Betti numbers are then
computed from the linear-dependence relations among the coefficient vectors of
the factors $z_0+\ell_\alpha$.

\subsection{Hyperplane Arrangements}
A complex hyperplane arrangement is a finite collection $\mathcal A=\{H_1,\ldots,H_m\}$ of affine hyperplanes in $\mathbb{C}^n$. Its complement is
\[
M(\mathcal A):=\mathbb{C}^n\setminus\bigcup_{H\in\mathcal A}H.
\]
The arrangement is central if every affine hyperplane contains the origin. Here, we will only be concerned with central arrangements. The Poincar\'{e} polynomial of the arrangement is
\[
P_{\mathcal A}(t):=\sum_{j\geq0}
\dim_{\mathbb{C}}H^j(M(\mathcal A),\mathbb{C})t^j.
\]
We associate with each hyperplane $H_j$ a symbol $e_j$ and let $E$ be the
exterior algebra generated by these symbols; thus
$e_i\wedge e_j=-e_j\wedge e_i$ and $e_i\wedge e_i=0$. For a defining linear form
$\ell_j$ with $H_j=\{\ell_j=0\}$, send $e_j$ to the de Rham cohomology class
of $(2\pi i)^{-1}d\ell_j/\ell_j$. This gives a graded algebra homomorphism
$\tau:E\to H^*(M(\mathcal A),\mathbb{C})$. Define the boundary map
$\partial:E\to E$ on a monomial by
\[
\partial(e_{i_1}\wedge\cdots\wedge e_{i_r})
=\sum_{j=1}^r(-1)^{j-1}e_{i_1}\wedge\cdots\wedge
\widehat{e_{i_j}}\wedge\cdots\wedge e_{i_r}.
\]
The Orlik--Solomon theorem identifies the kernel of $\tau$. Let $J$ be the
ideal generated by $\partial(e_{i_1}\wedge\cdots\wedge e_{i_r})$ for those
sets $\{\ell_{i_1},\ldots,\ell_{i_r}\}$ that are linearly dependent. The Orlik--Solomon
algebra is the quotient $A(\mathcal A):=E/J$. The following theorem lays the groundwork for the discussion in this section.

\begin{theorem}[\cite{OrlikSolomon1980}, Theorem 5.2]\label{thm:orlik-solomon}
There is an isomorphism of graded algebras $A(\mathcal A)\cong H^*(M(\mathcal A),\mathbb{C})$ carried by the map $\tau$.
\end{theorem}

This result equips us with a tool to study the topological consequences of
Theorem~\ref{thm:nilradical-factorization}. For more information about hyperplane arrangements, we refer the reader to Arnold~\cite{Arnold1969}, Brieskorn~\cite{Brieskorn1973}, and
Orlik--Solomon~\cite{OrlikSolomon1980}. 

Assume $\mathfrak{s}=\mathfrak{n}\oplus \mathfrak{v}$ is a solvable Lie algebra with $\dim \mathfrak{v}=k$. Write $Q_{\mathfrak{s}}=z_0^a p_{\mathfrak{s}}$, where $a=e_0(\mathfrak{s})$ and $z_0\nmid p_{\mathfrak{s}}$. Theorem \ref{thm:nilradical-factorization} shows that $p_{\mathfrak{s}}$ is a product of linear factors of the form $L_\alpha:=z_0+\ell_\alpha(\mathbf{z}_{\mathfrak v})$, $\alpha\in\Delta_{\mathfrak s}^\times$. We define $\Phi_{\mathfrak{s}}:=\{L_\alpha:\alpha\in\Delta_{\mathfrak s}^\times\}\subset (\mathbb{C}\oplus \mathfrak{v})^*$, $\mathcal A_{\mathfrak{s}}:=\{\ker L:L\in\Phi_{\mathfrak{s}}\}$, and set
\[M(\mathcal A_{\mathfrak{s}})=
\mathbb{C}^{k+1}\setminus\bigcup_{H\in\mathcal A_{\mathfrak{s}}}H.
\]
Thus $P_{\mathcal A_{\mathfrak s}}(t)$ denotes the Poincar\'{e} polynomial
of this arrangement complement.

Let $E_{\mathfrak s}$ be the exterior algebra generated by $\{e_L\}_{L\in\Phi_{\mathfrak{s}}}$ and $J_{\mathfrak s}$ be the Orlik--Solomon ideal of the arrangement $\mathcal A_{\mathfrak{s}}$. Then Theorem \ref{thm:orlik-solomon} determines the cohomology of the arrangement complement $M(\mathcal A_{\mathfrak{s}})$; $E_{\mathfrak s}/J_{\mathfrak s}$ is the Orlik--Solomon algebra of
$\mathcal A_{\mathfrak{s}}$ by definition, so Theorem~\ref{thm:orlik-solomon}
identifies it with the cohomology ring of the complement.

\begin{lemma}\label{lem:eigenvariety-product}
Let $\mathfrak{s}=\mathfrak{n}\oplus\mathfrak{v}$ be a solvable Lie algebra
with $k=\dim\mathfrak v$.  Under the coordinate decomposition
\[
\mathbb{C}^{\dim\mathfrak s+1}
\cong \mathbb{C}^{\dim\mathfrak n}\times\mathbb{C}^{k+1},
\]
the eigenvariety complement satisfies
\[
(V_{\mathfrak s})^c
\cong \mathbb{C}^{\dim\mathfrak n}\times M(\mathcal A_{\mathfrak s}).
\]
Consequently, projection onto the second factor is a homotopy equivalence
and induces an isomorphism
\[
H^*((V_{\mathfrak s})^c,\mathbb{C})
\cong H^*(M(\mathcal A_{\mathfrak s}),\mathbb{C}).
\]
\end{lemma}

\begin{proof}
By Theorem~\ref{thm:nilradical-factorization}, after the factor $z_0$ has
been removed, the polynomial $p_{\mathfrak s}$ depends only on $z_0$ and
the complement coordinates $\mathbf{z}_{\mathfrak v}$.  It is independent of
the $\dim\mathfrak n$ coordinates corresponding to $\mathfrak n$.  Hence the
complement splits as the displayed product.  The first factor is contractible,
so projection onto the second factor is a homotopy equivalence.
\end{proof}

The lemma identifies the cohomology of the eigenvariety complement
with that of $M(\mathcal A_{\mathfrak s})$.  By
Theorem~\ref{thm:orlik-solomon}, the latter is isomorphic to
$E_{\mathfrak s}/J_{\mathfrak s}$.  We record the resulting identification as a
corollary.

\begin{corollary}\label{cor:spectral-os}
Let $\mathfrak{s}=\mathfrak{n}\oplus \mathfrak{v}$ be a solvable Lie algebra.
Then
\[
H^*((V_{\mathfrak s})^c,\mathbb{C})
\cong H^*(M(\mathcal A_{\mathfrak s}),\mathbb{C})
\cong E_{\mathfrak s}/J_{\mathfrak s}.
\]
\end{corollary}
\subsection{Spectral Matroid and the Poincar\'{e} Polynomial}
Corollary \ref{cor:spectral-os} reduces the computation of the cohomology ring to the combinatorics of linear dependencies, which can be studied through the theory of matroids. For any linear function $z_0+c_1z_{f_1}+\cdots +c_kz_{f_k}$, its coefficient vector is $(1,c_1,\ldots,c_k)\in \mathbb{C}^{k+1}$. We begin by recalling some basics about matroids. For a more comprehensive introduction, see \cite{Oxley2011}.

\begin{definition}
A \emph{matroid} is a pair $\mathcal M=(E,\mathcal I)$ consisting of a finite
ground set $E$ and a collection $\mathcal I$ of subsets of $E$, called the
independent sets, such that:
\begin{enumerate}
\item $\varnothing\in\mathcal I$;
\item if $I\in\mathcal I$ and $J\subseteq I$, then $J\in\mathcal I$;
\item if $I,J\in\mathcal I$ and $|I|<|J|$, then there is an element
$e\in J\setminus I$ for which $I\cup\{e\}\in\mathcal I$.
\end{enumerate}
\end{definition}

\begin{definition}
    For a matroid $\mathcal M=(E,\mathcal I)$ and a subset $A\subseteq E$, its
rank is
\[
r_{\mathcal M}(A):=\max\{|I|: I\subseteq A,\ I\in\mathcal I\}.
\]
The rank of $\mathcal M$ is $r_{\mathcal M}(E)$, denoted
$\operatorname{rank}\mathcal M$.
\end{definition}

A {\em circuit} is a minimal dependent subset of the ground set. An element is a
{\em coloop} if it belongs to no circuit. Here a subset is dependent when it is not
independent.

\begin{definition}
Two matroids $\mathcal M=(E,\mathcal I)$ and
$\mathcal M'=(E',\mathcal I')$ are said to be isomorphic, denoted $\mathcal M\cong\mathcal M'$, if there exists a bijection
$\varphi:E\to E'$ such that $I\in\mathcal I$ if and only if
$\varphi(I)\in\mathcal I'$ for every $I\subseteq E$.
\end{definition}

Having recalled the necessary terminology, we now associate a matroid to our setting.

\begin{definition}
Let $\mathfrak{s}=\mathfrak{n}\oplus \mathfrak{v}$ be a solvable Lie algebra.
The \emph{spectral matroid} is
$\mathcal M_{\mathfrak s}=(\Phi_{\mathfrak s},\mathcal I_{\mathfrak s})$,
where the ground set $\Phi_{\mathfrak s}$ consists of the distinct non-$z_0$
linear factors and $\mathcal I_{\mathfrak s}$ consists of the subsets whose
coefficient vectors are linearly independent.
\end{definition}

This matroid is associated with the Lie algebra through its characteristic
polynomial, not to an arbitrary choice of arrangement. The factors in
$\Phi_{\mathfrak{s}}$ are the distinct factors other than $z_0$ supplied by
Theorem~\ref{thm:nilradical-factorization}; changing the adapted basis applies
an invertible linear transformation to their coefficient vectors and therefore
does not change the represented matroid. The multiplicities remain part of
the characteristic polynomial, while the matroid records the additional
information carried by the linear-dependence relations among their coefficient
vectors. Equivalently, these are the affine-dependence relations
among the distinct nonzero nilradical weights. Thus
$\mathcal M_{\mathfrak{s}}$ is a coordinate-independent record of this
dependence structure. The following lemma shows that it is also invariant under
spectral equivalence.

\begin{lemma}\label{lem:spectral-matroid-iso}
Suppose $\mathfrak{s}$ and $\mathfrak{s}'$ are solvable complex Lie algebras. If $\mathfrak{s}\sim\mathfrak{s}'$, then $\mathcal M_{\mathfrak{s}}\cong\mathcal M_{\mathfrak{s}'}$.
\end{lemma}

\begin{proof}
 Suppose $\mathfrak{s}\sim\mathfrak{s}'$. By Proposition~\ref{prop:weight-equivalence}, the total $z_0$-multiplicities agree and there is a linear isomorphism $T$ between the spans of the nonzero nilradical weight multisets carrying $\Delta_{\mathfrak{s}}^\times$ to $\Delta_{\mathfrak{s}'}^\times$ as multisets.

The ground set $\Phi_{\mathfrak{s}}$ consists of the linear functions $z_0+\ell_\alpha$, where $\alpha$ ranges over the distinct nonzero nilradical weights. The map $z_0+\ell_\alpha\mapsto z_0+\ell_{T(\alpha)}$ gives a bijection $\Phi_{\mathfrak{s}}\to\Phi_{\mathfrak{s}'}$. A subset $\{z_0+\ell_{\alpha_1},\dots,z_0+\ell_{\alpha_r}\}$ is linearly dependent if and only if there are coefficients $c_1,\dots,c_r$, not all zero, such that $\sum_i c_i=0$ and $\sum_i c_i\alpha_i=0$. Since $T$ is a linear isomorphism on the span of the weights, this condition is equivalent to $\sum_i c_i=0$ and $\sum_i c_iT(\alpha_i)=0$. Thus the bijection preserves the dependent subsets and induces an isomorphism $\mathcal M_{\mathfrak{s}}\cong\mathcal M_{\mathfrak{s}'}$.
\end{proof}

To compute the Betti numbers of $V_{\mathfrak{s}}^c$ explicitly, we introduce the tools related to the intersection lattice.

\begin{definition}
Let $\mathcal A_{\mathfrak{s}}\subset \mathbb{C}^{k+1}$ be the hyperplane arrangement associated with a solvable Lie algebra $\mathfrak{s}$. The \emph{intersection lattice} $L(\mathcal A_{\mathfrak{s}})$ consists of all subspaces $X\subseteq \mathbb{C}^{k+1}$ of the form
\[
X=\bigcap_{H\in S}H
\]
for some $S\subseteq \mathcal A_{\mathfrak{s}}$. The partial order is given by reverse inclusion: $X\leq Y$ if $Y\subseteq X$. The minimal element is $\hat 0=\mathbb{C}^{k+1}$, and the maximal element is $\hat 1=\bigcap_{H\in\mathcal A_{\mathfrak{s}}}H$. The rank function is defined by $\operatorname{rank}(X)=\operatorname{codim}_{\mathbb{C}^{k+1}} X$, where codimension is the difference between the ambient dimension and $\dim X$. For the empty arrangement, we use the convention
$\bigcap_{H\in\varnothing}H=\mathbb{C}^{k+1}$, so the intersection lattice has the single element $\hat 0=\hat 1=\mathbb{C}^{k+1}$.
\end{definition}

The inclusion--exclusion structure of the arrangement is given by the Möbius function of this lattice; see, for example, \cite{Oxley2011}.

\begin{definition}
    The \emph{Möbius function} $\mu \colon L(\mathcal{A}_{\mathfrak{s}}) \times L(\mathcal{A}_{\mathfrak{s}}) \to \mathbb{Z}$ is characterized by $\mu(X,X)=1$ for all $X$, and by the relation
    \[
        \sum_{X \le Z \le Y} \mu(X,Z) = 0 \quad \text{for } X<Y.
    \]
    We write $\mu(X)=\mu(\hat{0},X)$.
\end{definition}

To express the Betti numbers combinatorially, we introduce the Whitney numbers
of the first kind; see \cite{Oxley2011}.

\begin{definition}
    The \emph{Whitney numbers of the first kind} are defined by
    \[
        w_k = \sum_{\operatorname{rank}(X)=k} |\mu(X)|.
    \]
\end{definition}

Since $\mathcal A_{\mathfrak{s}}$ is a central complex hyperplane arrangement,
its intersection lattice is a finite geometric lattice, meaning the lattice of
flats associated with a matroid; in the represented matroid here, a flat is a
set of coefficient vectors containing every ground-set vector in its span.
We use the following standard form of the
Orlik--Solomon formula.

\begin{theorem}[\cite{OrlikSolomon1980}, Theorem 2.6]\label{thm:os-betti}
Let $L$ be a finite geometric lattice. The Poincar\'{e} polynomial of the associated
graded Orlik--Solomon algebra is
\[
P(t)=\sum_{X\in L}\mu(X)(-t)^{\operatorname{rank}(X)},
\]
where $\mu(X)=\mu(\hat 0,X)$.
\end{theorem}

Let $w_k$ be the Whitney numbers of the first kind of the intersection lattice
$L(\mathcal A_{\mathfrak{s}})$. The theorem above gives the Poincar\'{e} polynomial
of the arrangement complement $M(\mathcal A_{\mathfrak{s}})$.

\begin{corollary}\label{cor:poincare-matroid}
Let $\mathcal A_{\mathfrak{s}}$ be the central arrangement associated with a
solvable Lie algebra $\mathfrak{s}$.  Then the Poincar\'{e} polynomial of
$V_{\mathfrak{s}}^c$ (equivalently, of $M(\mathcal A_{\mathfrak{s}})$) is
\[
P_{\mathfrak{s}}(t)=\sum_{k\geq0}w_kt^k.
\]
Equivalently, the Betti numbers $b_k(\mathfrak s)=w_k$ for every $k\geq0$.
\end{corollary}

\begin{proof}
By Lemma~\ref{lem:eigenvariety-product} and Corollary~\ref{cor:spectral-os},
the cohomology of $V_{\mathfrak{s}}^c$ is isomorphic to the
Orlik--Solomon algebra $A(\mathcal A_{\mathfrak{s}})$.  Theorem
\ref{thm:os-betti} gives
\[
P_{\mathfrak{s}}(t)
=
\sum_{X\in L(\mathcal A_{\mathfrak{s}})}
\mu(X)(-t)^{\operatorname{rank}(X)}.
\]
For a geometric lattice, $\mu(X)$ has sign
$(-1)^{\operatorname{rank}(X)}$. Therefore, the coefficient of $t^k$ is
$\sum_{\operatorname{rank}(X)=k}|\mu(X)|=w_k$.
\end{proof}

Thus the Betti numbers $b_j(\mathfrak{s})$ are determined by the spectral
matroid $\mathcal{M}_{\mathfrak{s}}$, providing an answer to Problem~\ref{prob:higher-betti}.

We may now determine the degree of $P_{\mathfrak{s}}(t)$ from the nilradical weights. This sharpens the earlier bound in Proposition~\ref{prop:first-betti} by identifying the rank that gives the degree. Recall that $\Delta_{\mathfrak s}^{\times}$ denotes the nonzero nilradical weight multiset and $\Phi_{\mathfrak s}$ is the set of distinct associated factors $L_\alpha=z_0+\ell_\alpha$. For a nonempty finite set $B=\{{\bf b}_1,\ldots,{\bf b}_p\}$ of weights, define its affine span by
\[
\operatorname{Aff}(B)
:=
\left\{\sum_{i=1}^p t_i b_i:
t_i\in\mathbb{C},\ \sum_{i=1}^p t_i=1\right\}.
\]

\begin{theorem}\label{thm:poincare-degree}
    If $\Delta_{\mathfrak s}^{\times}=\varnothing$, then $p_{\mathfrak s}=1$ and
    $P_{\mathfrak s}(t)=1$. Otherwise,
    \[
    \deg P_{\mathfrak s}(t)
    =
    \operatorname{rank}\mathcal M_{\mathfrak s}
    =
    \dim\operatorname{span}\Phi_{\mathfrak s}
    =
    \operatorname{codim}
    \bigcap_{L_\alpha\in\Phi_{\mathfrak s}}\ker L_\alpha.
    \]
    Moreover, for any $\alpha_0\in\Delta_{\mathfrak s}^{\times}$,
    \[
    \deg P_{\mathfrak s}(t)
    =
    1+\dim\operatorname{span}
    \{\alpha-\alpha_0:\alpha\in\Delta_{\mathfrak s}^{\times}\}
    =
    1+\dim\operatorname{Aff}(\Delta_{\mathfrak s}^{\times}).
    \]
    \end{theorem}

    \begin{proof}
    By Corollary~\ref{cor:poincare-matroid},
    $P_{\mathfrak s}(t)$ is the Poincar\'e polynomial of the
    Orlik--Solomon algebra of the matroid represented by
    $\Phi_{\mathfrak s}$. Its degree is therefore the rank of this matroid, namely
    $\dim\operatorname{span}\Phi_{\mathfrak s}$. Since the forms in
    $\Phi_{\mathfrak s}$ define the associated arrangement, this dimension is also
    \[
    \operatorname{codim}
    \bigcap_{L_\alpha\in\Phi_{\mathfrak s}}\ker L_\alpha.
    \]

    If $\Delta_{\mathfrak s}^{\times}=\varnothing$, there are no non-$z_0$ factors, so
    $p_{\mathfrak s}=1$ and $P_{\mathfrak s}(t)=1$.

    Now suppose $\Delta_{\mathfrak s}^{\times}\neq\varnothing$ and fix
    $\alpha_0\in\Delta_{\mathfrak s}^{\times}$. Since
    \[
    L_\alpha=L_{\alpha_0}+\ell_{\alpha-\alpha_0},
    \]
    we have
    \[
    \operatorname{span}\Phi_{\mathfrak s}
    =
    \mathbb{C} L_{\alpha_0}
    \oplus
    \operatorname{span}
    \{\ell_{\alpha-\alpha_0}:\alpha\in\Delta_{\mathfrak s}^{\times}\}.
    \]
    The sum is direct because the second summand lies in $\mathfrak v^*$, while
    $L_{\alpha_0}$ has nonzero $z_0$-component. Taking dimensions gives the result.
    \end{proof}

\begin{remark}
It is important that the degree depends on the affine span of the nonzero weights, not on their ordinary linear span. For example, if $\dim \mathfrak{v}=1$ and $\Delta_{\mathfrak s}^\times=\{\alpha\}$ with $\alpha\neq0$, then $\Delta_{\mathfrak s}^\times$ spans $\mathfrak{v}^*$ as a vector space, but $\operatorname{rank}\{L_\alpha\}=1$, so $P_{\mathfrak s}(t)$ has degree $1$, not $2$.
\end{remark}

\subsection{The Log-Concavity of the Betti Numbers}
Theorem~\ref{thm:poincare-degree} gives a direct upper bound for the Betti numbers in terms of the spectral index $\kappa(\mathfrak{s})$, equivalently the number of elements in the ground set $\Phi_{\mathfrak{s}}$ of the spectral matroid.
\begin{corollary}\label{cor:betti-bound}
For any solvable Lie algebra $\mathfrak{s}$ and $j\geq0$, we have
\[
b_j(\mathfrak{s})\leq \binom{\kappa(\mathfrak{s})-1}{j}.
\]
\end{corollary}

\begin{proof}
By Corollary~\ref{cor:spectral-os}, the cohomology
of $V_{\mathfrak{s}}^c$ is isomorphic to the Orlik--Solomon algebra $E_\mathfrak{s}/J_\mathfrak{s}$. Therefore, the degree-$j$ cohomology space $H^j(V_{\mathfrak{s}}^c,\mathbb{C})$ is a quotient of the degree-$j$ part of $E_\mathfrak{s}$, which has dimension
$\binom{\kappa(\mathfrak{s})-1}{j}$. The corollary follows.
\end{proof}

We proceed by establishing a qualitative property of the Betti
numbers of $V_{\mathfrak{s}}^c$. Recall that a sequence of non-negative real numbers $a_0,a_1,\ldots,a_n$ is log-concave if $a_k^2\geq a_{k-1}a_{k+1}$ for $1\leq k<n$. The sequence is said to have \emph{no internal zeros} if every entry between two nonzero entries is nonzero, and it is \emph{unimodal} if it weakly increases up to some index and then weakly decreases. If a log-concave sequence has no internal zeros, then it is unimodal. To prove log-concavity
for the Betti numbers, we introduce the characteristic polynomial of a spectral
matroid.

\begin{definition}
    Let $\mathcal{M}$ be a matroid of rank $r$. The \emph{characteristic polynomial} of $\mathcal{M}$ is defined via the Möbius function of its lattice as
    \[ \chi_{\mathcal{M}}(t) = \sum_{X \in L(\mathcal{M})} \mu(\hat{0}, X) t^{r - \operatorname{rank}(X)}. \]
\end{definition}

The following result concerns the coefficients of $\chi_{\mathcal{M}}$.

\begin{theorem}[\cite{HuhKatz2012}]\label{thm:huh-katz}
Let $\mathcal M$ be a matroid representable over a field of characteristic $0$. Then the absolute values of the coefficients of its characteristic polynomial form a log-concave sequence.
\end{theorem}

Our next lemma connects the Poincar\'{e} polynomial of a hyperplane arrangement with the characteristic polynomial of a matroid.

\begin{lemma}\label{lem:poly-reversal}
Let $\mathfrak{s}$ be a solvable Lie algebra with $\operatorname{rank}\mathcal M_{\mathfrak{s}}=r$. Then \[\chi_{\mathcal M_{\mathfrak{s}}}(q)=q^rP_{\mathfrak{s}}(-q^{-1}).\]
\end{lemma}

\begin{proof}
Using the definition of the Whitney numbers
\[
w_k=\sum_{\operatorname{rank}(X)=k}|\mu(\hat 0,X)|,
\]
and the sign property $\mu(\hat 0,X)=(-1)^{\operatorname{rank}(X)}|\mu(\hat 0,X)|$ for geometric lattices, the characteristic polynomial of $\mathcal M_{\mathfrak{s}}$ is
\[
\chi_{\mathcal M_{\mathfrak{s}}}(q)
=
\sum_{k=0}^r (-1)^k w_k q^{r-k}.
\]
Since
\[
P_{\mathfrak{s}}(t)=\sum_{k=0}^r w_kt^k,
\]
we obtain
\[
q^rP_{\mathfrak{s}}(-q^{-1})
=
q^r\sum_{k=0}^r w_k(-q^{-1})^k
=
\sum_{k=0}^r (-1)^k w_k q^{r-k}
=
\chi_{\mathcal M_{\mathfrak{s}}}(q).
\]
\end{proof}
Therefore, the Betti numbers of $V_{\mathfrak{s}}^c$ coincide with the absolute values of the coefficients of the characteristic polynomial $\chi_{\mathcal M_{\mathfrak{s}}}(q)$. The following is a direct consequence of Theorem~\ref{thm:huh-katz}.

\begin{corollary}\label{cor:betti-log-concave}
Let $\mathfrak{s}$ be a solvable Lie algebra with $\operatorname{rank}\mathcal M_{\mathfrak{s}}=r$. Then the Betti numbers $b_0(\mathfrak{s}),\ldots,b_r(\mathfrak{s})$ form a log-concave sequence:
\[
b_{j-1}(\mathfrak{s})\,b_{j+1}(\mathfrak{s})\leq b_j(\mathfrak{s})^2,
\qquad 1\leq j\leq r-1.
\]
Moreover, $b_j(\mathfrak{s})=0$ for $j>r$.
\end{corollary}

\begin{proof}
The spectral matroid $\mathcal M_{\mathfrak{s}}$ is represented by linear
forms over $\mathbb{C}$, hence is representable over a field of characteristic
$0$. By Theorem~\ref{thm:huh-katz}, the absolute values of the coefficients of
$\chi_{\mathcal M_{\mathfrak{s}}}(q)$ form a log-concave sequence. By
Lemma~\ref{lem:poly-reversal}, these absolute values are the Betti numbers
$b_j(\mathfrak{s})$ for $0\leq j\leq r$. The vanishing for $j>r$ follows from
Theorem~\ref{thm:poincare-degree}.
\end{proof}

\section{Solvable Lie Algebras with Abelian Nilradical}
\label{sec:abelian-nilradical}

Assume $\mathfrak{s}$ is a solvable Lie algebra with an abelian nilradical $\mathfrak{n}$ and suppose the nilradical weight multiset $\Delta_{\mathfrak{s}}$ consists of the weights $\mathbf a_j$, each repeated $r_j$ times. The form of such $\mathfrak{s}$ is described by Ndogmo and Winternitz \cite{NdogmoWinternitz1994}. The calculations
below isolate the weights contributed by the nilradical blocks from
the remaining extension data.

\subsection{Classification Theorems of Ndogmo and Winternitz}

Write $r=\dim\mathfrak n$ and $k=\dim(\mathfrak s/\mathfrak n)$, so
$\dim\mathfrak s=r+k$ and $1\leq k\leq r$.  Ndogmo and Winternitz \cite{NdogmoWinternitz1994} reduce the
isomorphism problem to a commuting family of matrices together with compatible
extension data, represented here by the brackets between a chosen complement
and itself.  One must then classify this extension data modulo the relevant
changes of basis and remove decomposable algebras. Here, decomposable means
isomorphic to a direct sum of two nonzero ideals.
Choose a vector-space complement $\mathfrak v$ to $\mathfrak n$ and bases
$\{n_1,\ldots,n_r\}$ of $\mathfrak n$ and $\{x_1,\ldots,x_k\}$ of
$\mathfrak v$.  Thus $\mathfrak s=\mathfrak n\oplus\mathfrak v$ as a vector
space.  Let $z_1,\ldots,z_r$ be the variables corresponding to
$(n_1,\ldots,n_r)$, let $z_{r+\alpha}$ correspond to $x_\alpha$, and write
$\mathbf{z}_{\mathfrak v}=(z_{r+1},\ldots,z_{r+k})$.
The point of the calculation below is that the variables corresponding to
$n_1,\ldots,n_r$, the extension data, and the nilpotent parts of the
commuting matrices all disappear from the determinant.  Thus the determinant
depends only on the joint eigenvalue data of the commuting family, equivalently
on the simultaneous diagonal weights and their multiplicities.

For $1\leq\alpha\leq k$, define $A^\alpha\in\mathfrak{gl}(r,\mathbb{C})$ and
$R_{\alpha\beta}\in\mathfrak n$ by
$[x_\alpha,n_j]=\sum_{i=1}^r(A^\alpha)_{ij}n_i$ for $1\leq j\leq r$ and
$[x_\alpha,x_\beta]=R_{\alpha\beta}$.
The vectors $R_{\alpha\beta}$ measure the failure of the chosen complement
to be a Lie subalgebra, and their compatibility relations are the
additional conditions imposed by the Jacobi identity.

We call the commuting family $A^1,\ldots,A^k$ linearly nilindependent if
no nonzero linear combination $\sum_{\alpha=1}^k c_\alpha A^\alpha$ is
nilpotent.

\begin{theorem}
\label{thm:nw-structural-form}
Let $\mathfrak s$ be a finite-dimensional solvable non-nilpotent complex Lie
algebra with abelian nilradical $\mathfrak n$.  With respect to the bases
chosen above, the brackets have the form
\[
[x_\alpha,n_j]=\sum_{i=1}^r(A^\alpha)_{ij}n_i,\qquad
[n_i,n_j]=0,\qquad
[x_\alpha,x_\beta]=R_{\alpha\beta}\in\mathfrak n.
\]
The matrices $A^1,\ldots,A^k$ commute and are linearly nilindependent.  For
$k\geq3$ and
$1\leq\alpha<\beta<\gamma\leq k$,
\[
A^\alpha R_{\beta\gamma}+A^\beta R_{\gamma\alpha}
 +A^\gamma R_{\alpha\beta}=0.
\]
The classification is taken modulo the following transformations:
\begin{enumerate}
\item $x_\alpha\mapsto x_\alpha+t_\alpha$, with
      $t_\alpha\in\mathfrak n$;
\item $(x_1,\ldots,x_k)^{\mathsf T}\mapsto
      G(x_1,\ldots,x_k)^{\mathsf T}$, with $G\in\operatorname{GL}(k,\mathbb{C})$;
\item $(n_1,\ldots,n_r)^{\mathsf T}\mapsto
      S(n_1,\ldots,n_r)^{\mathsf T}$, with $S\in\operatorname{GL}(r,\mathbb{C})$.
\end{enumerate}
\end{theorem}

Thus the classification problem separates into the commuting matrix family and
the extension data $R_{\alpha\beta}$, followed by the removal of decomposable
algebras.  Over $\mathbb{C}$, the commuting family has the
following block form.

\begin{theorem}
\label{thm:nw-complex-block-form}
Let $\mathfrak s$ be as above over $\mathbb{C}$.  The matrices $A^\alpha$ can be
simultaneously transformed into the block diagonal form
\[
 A^\alpha=\bigoplus_{j=1}^{p+q}T_j^\alpha(a_j^\alpha),
 \qquad
 T_j^\alpha(a_j^\alpha)=a_j^\alpha I_{r_j}+N_j^\alpha,
 \qquad
 \sum_{j=1}^{p+q}r_j=r,
\]
where each $T_j^\alpha$ is triangular and the matrices $N_j^\alpha$ in a
fixed block span an abelian nilpotent subalgebra of
$\mathfrak{sl}(r_j,\mathbb{C})$.  Writing
$\mathbf a_j=(a_j^1,\ldots,a_j^k)$, the weights satisfy
\[
\mathbf a_j\neq0\quad(1\leq j\leq p),\qquad
\mathbf a_j=0\quad(p+1\leq j\leq p+q),
\]
and, after a change of basis in $\mathfrak s/\mathfrak n$,
\[
a_j^\alpha=\delta_j^\alpha\quad(1\leq j\leq k),
\qquad 1\leq k\leq p,\qquad p+q\leq r.
\]
\end{theorem}

This block decomposition reduces the classification of the commuting action to
the choice of a partition, the relevant indecomposable nilpotent block
subalgebras, and the compatible scalar and nilpotent parameter maps, modulo
changes of nilradical basis and complement basis.  Classification of the Lie
algebra additionally requires the compatible extension tensors
$R_{\alpha\beta}$ and the removal of decomposable cases.

The weights and block sizes enter the calculation of the characteristic
polynomial and the associated invariants.
The strictly triangular entries and the extension data remain important for
isomorphism, but not for the factorization considered here. These weights are the data used in the spectral classification below.

\subsection{Characteristic Polynomials in the Abelian Case}

With the notation of Theorem~\ref{thm:nw-complex-block-form}, write
$\mathbf{z}_{\mathfrak v}=(z_{r+1},\ldots,z_{r+k})$ for the complement
coordinates.
\begin{theorem}
\label{thm:abelian-block-factorization}
The characteristic polynomial of a solvable Lie algebra $\mathfrak{s}$ with
abelian nilradical is
\[
Q_{\mathfrak s}(z_0,z_1,\ldots,z_{r+k})
 =z_0^{e_0(\mathfrak s)}
   \prod_{j=1}^{p}
   \left(z_0+\sum_{\alpha=1}^k a_j^\alpha z_{r+\alpha}\right)^{r_j}.
\]
Moreover,
\[
e_0(\mathfrak s)=k+\sum_{j=p+1}^{p+q}r_j,
\qquad
\kappa(\mathfrak s)=1+|B|,
\]
where $B$ is the set of distinct nonzero weights among
$\mathbf a_1,\ldots,\mathbf a_p$.
\end{theorem}

\begin{proof}
Relative to $\mathfrak s=\mathfrak n\oplus\mathfrak v$, the matrix pencil in
the definition of $Q_{\mathfrak s}$ has the block form
\[
z_0I+\sum_{i=1}^r z_i\operatorname{ad}(n_i)
 +\sum_{\alpha=1}^k z_{r+\alpha}\operatorname{ad}(x_\alpha)
=
\begin{pmatrix}
\displaystyle z_0I_r+\sum_{\alpha=1}^k z_{r+\alpha} A^\alpha&*\\
0&z_0I_k
\end{pmatrix}.
\]
The lower-right block is $z_0I_k$ because
$[\mathfrak s,\mathfrak s]\subseteq\mathfrak n$, so the induced bracket on
$\mathfrak s/\mathfrak n$ vanishes.  Taking the determinant gives
\[
Q_{\mathfrak s}(z_0,z_1,\ldots,z_{r+k})
=z_0^k\det\left(z_0I_r+\sum_{\alpha=1}^k z_{r+\alpha} A^\alpha\right).
\]

On the $j$th block, the diagonal of
$\sum_\alpha z_{r+\alpha}A^\alpha$ is the scalar
$\sum_\alpha a_j^\alpha z_{r+\alpha}$, repeated $r_j$ times.
The triangular form therefore contributes
\[
\left(z_0+\sum_{\alpha=1}^k a_j^\alpha z_{r+\alpha}\right)^{r_j}.
\]
For $j>p$ the weight is zero, so the contribution is $z_0^{r_j}$.
Multiplying the quotient contribution with the contributions of all blocks
proves the factorization.  The factor $z_0$ contributed by the quotient block
is always present, and zero-weight blocks contribute additional copies of the
same factor.  Thus the distinct linear factors are indexed by $0$ together
with the distinct nonzero weights in $B$, which gives
$\kappa(\mathfrak s)=1+|B|$.
\end{proof}

\subsection{Classification up to Spectral Equivalence}

Theorem~\ref{thm:abelian-block-factorization} reduces spectral equivalence
to these weights.  The spectral equivalence class is specified by the
multiplicity of the $z_0$ factor and the multiset of nonzero weights with
their multiplicities, modulo a change of basis in the extension space.  The
remaining normal-form parameters continue to matter for isomorphism.
By Proposition~\ref{prop:weight-equivalence}, this description is equivalent
to equality of the $z_0$-multiplicities together with linear equivalence of
the nonzero weight multisets.  We record the resulting specialization.

\begin{remark}
\label{cor:nw-spectral-equivalence}
In the abelian-nilradical setting, Proposition~\ref{prop:weight-equivalence}
applies directly: spectral equivalence is controlled by $e_0(\mathfrak s)$ and
the linear-equivalence class of the nonzero nilradical weight multiset.
\end{remark}

The criterion is also sharp: every finite multiset of nonzero weights whose
span is $\mathbb{C}^k$, together with any prescribed nonnegative zero-weight
multiplicity, is realizable.

\begin{proposition}
\label{prop:abelian-weight-realization}
Assume $B:=\{\mathbf b_1,\ldots,\mathbf b_p\}\subset \mathbb{C}^k\setminus\{0\}$ and it spans
$\mathbb{C}^k$. Let $r_0,\ldots,r_p$ be nonnegative integers with $r_j>0$ for each $j\geq 1$
and set $r=r_0+\cdots +r_p$.
Then there is a solvable complex Lie algebra $\mathfrak{s}$ with abelian nilradical $\mathfrak n$ and
$\dim(\mathfrak s/\mathfrak n)=k$ for which
\[
Q_{\mathfrak s}(z_0,z_1,\ldots,z_{r+k})
 =z_0^{k+r_0}\prod_{i=1}^p
 \left(z_0+\sum_{\alpha=1}^k b_i^\alpha z_{r+\alpha}\right)^{r_i}.
\]
\end{proposition}

\begin{proof}
Let $\mathfrak n=N_0\oplus N_1\oplus\cdots\oplus N_p$, where
$\dim N_i=r_i$, and let
$\mathfrak v=\operatorname{span}\{x_1,\ldots,x_k\}$. Define
$A^\alpha|_{N_0}=0$ and
$A^\alpha|_{N_i}=b_i^\alpha I_{N_i}$, and define
$\mathfrak s=\mathfrak n\oplus\mathfrak v$ by
$[x_\alpha,u]=A^\alpha u$ for $u\in\mathfrak n$ and
$[x_\alpha,x_\beta]=0$. The matrices are diagonal with respect to this
decomposition, so they commute, and the choice $R_{\alpha\beta}=0$ satisfies
the extension compatibility conditions. Thus the Jacobi identity holds. The
extension tensors do not enter the characteristic polynomial in
Theorem~\ref{thm:abelian-block-factorization}, although nonzero compatible
choices can affect the isomorphism class.

Suppose that $A(c):=\sum_\alpha c_\alpha A^\alpha$ is nilpotent.  On $N_i$ it
acts by the scalar $\sum_\alpha c_\alpha b_i^\alpha$.  A scalar operator is
nilpotent only when its scalar is zero; hence
$\sum_\alpha c_\alpha b_i^\alpha=0$ for every $i$.  Since the vectors
$\mathbf b_i$ span $\mathbb{C}^k$, this forces $c_1=\cdots=c_k=0$.
Thus $A^1,\ldots,A^k$ are linearly nilindependent.  Let
$\widetilde{\mathfrak n}$ denote the nilradical of the constructed algebra.
Since $\mathfrak n$ is an abelian ideal,
$\mathfrak n\subseteq\widetilde{\mathfrak n}$.  If
$\widetilde{\mathfrak n}/\mathfrak n\neq0$, choose
\[
y=u+\sum_{\alpha=1}^k c_\alpha x_\alpha
\in\widetilde{\mathfrak n}\setminus\mathfrak n.
\]
Because $\widetilde{\mathfrak n}$ is a nilpotent ideal,
$[y,\mathfrak s]\subseteq\widetilde{\mathfrak n}$ and
$\operatorname{ad}(y)|_{\widetilde{\mathfrak n}}$ is nilpotent.  Thus
$\operatorname{ad}(y)$ induces the zero operator on
$\mathfrak s/\widetilde{\mathfrak n}$, since
$[\mathfrak s,\mathfrak s]\subseteq\mathfrak n\subseteq\widetilde{\mathfrak n}$.
Relative to $\widetilde{\mathfrak n}\subseteq\mathfrak s$, it is therefore
block upper triangular with nilpotent diagonal blocks, and is nilpotent on
$\mathfrak s$.  Hence its restriction
$A(c)=\sum_\alpha c_\alpha A^\alpha$ to $\mathfrak n$ is nilpotent.
Linear nilindependence forces $c_1=\cdots=c_k=0$, a contradiction.  Thus
$\widetilde{\mathfrak n}=\mathfrak n$, and Theorem~\ref{thm:abelian-block-factorization}
applied to these blocks gives the displayed characteristic polynomial.
\end{proof}

Note that we do not claim that the construction is indecomposable. In particular, if
$r_0>0$ and all $R_{\alpha\beta}$ vanish, then $N_0$ is a central direct
summand.

\subsection{Topological Invariants in the Abelian Case}

In light of Theorem \ref{thm:poincare-degree}, the degree of the Poincar\'{e} polynomial of the hyperplane complement $M({\mathcal A}_\mathfrak{s})$ is determined by ordinary affine dependence among the distinct nonzero weights. For $\mathbf b\in\mathbb{C}^k$, define
$\widehat{\mathbf b}:=(1,\mathbf b)$. For a finite set
$B\subset\mathbb{C}^k$ with $|B|=p$, we call a subset affinely independent when its lifted
vectors $\widehat{\mathbf b}$ are linearly independent. We say that $B$ is in
affine general position if every subset of at most $m=\min\{k+1,p\}$ points is affinely independent.

For integers $0\leq m\leq p$, let $U_{m,p}$ denote the uniform
matroid on a $p$-element ground set whose independent subsets are precisely
those of cardinality at most $m$.

\begin{corollary}
\label{cor:abelian-affine-matroid}
Let $\mathfrak{s}$ be a solvable Lie algebra determined by Proposition \ref{prop:abelian-weight-realization}. Then the spectral matroid is represented by the vectors
$\widehat{\mathbf b}=(1,b^1,\ldots,b^k)$ for $\mathbf b\in B$, and
\[
 \deg P_{\mathfrak s}(t)=1+\dim\operatorname{Aff}(B).
\]
If $B$ is in affine general position, then the spectral matroid is
$U_{m,p}$, and
\[
 P_{\mathfrak s}(t)=
 \sum_{i=0}^{m-1}\binom{p}{i}t^i
 +\binom{p-1}{m-1}t^m.
\]
\end{corollary}

\begin{proof}
The coefficient vector of
$z_0+\sum_\alpha b^\alpha z_{r+\alpha}$ is
$\widehat{\mathbf b}$.  Linear dependence among these coefficient vectors is
equivalent to affine dependence among the points of $B$.  Thus a subset of
the displayed linear forms is dependent if and only if the corresponding
weights are affinely dependent.  The rank and degree formula follows from
Theorem~\ref{thm:poincare-degree}.  In affine general position every subset of at
most $m$ weights is independent, so the represented matroid is the uniform matroid
$U_{m,p}$.  Its Whitney numbers are
$\binom{p}{i}$ for $i<m$ and $\binom{p-1}{m-1}$ for $i=m$.
\end{proof}

\begin{remark}
The affine dependence appears because the coefficient vector of the linear
form $z_0+\sum_\alpha b^\alpha z_{r+\alpha}$ is
$\widehat{\mathbf b}=(1,b^1,\ldots,b^k)$.  The first coordinate is the same
for every weight, so linear dependence among the coefficient vectors
$\widehat{\mathbf b}$ is equivalent to affine dependence among the weights
$\mathbf b$.  The arrangement complement therefore records the
affine-dependence structure of the weight set.
\end{remark}

The general bounds for the degree and coefficients of $P_{\mathfrak s}$ therefore specialize to
\[
\deg P_{\mathfrak s}\leq k+1,
\qquad
b_i(\mathfrak s)\leq\binom{|B|}{i}
=\binom{\kappa(\mathfrak s)-1}{i}.
\]

For a finite set $B$ of distinct nonzero weights, let
\[
\mathcal A_B:=\{\ker(z_0+\ell_{\mathbf b}):\mathbf b\in B\}.
\]
Thus
$P_{\mathcal A_B}(t)$ is the polynomial defined above for this arrangement;
the same notation applies to $P_{\mathcal A_S}(t)$ and
$P_{\mathcal A_{B_{\mathbf a}}}(t)$ below. This notation is used in
Tables~\ref{tab:nw-maximal-center-spectra} and~\ref{tab:nw-853-spectra}.

\subsection{Tables of Weights and Spectral Invariants}

We apply the preceding description to the special families studied in
\cite{NdogmoWinternitz1994}.  In each case, we record the weight multiset, the
resulting characteristic polynomial, and selected matroid invariants.

The tables record characteristic-polynomial formulas and selected matroid
invariants. Equal values of $\kappa$ and $P_{\mathfrak s}$ do not in general
imply spectral equivalence; the full linear-equivalence class of the weight
multiset is still required. For example, in the case $k=1$, the weight sets
$\{1,2,3\}$ and $\{1,2,4\}$ have the same uniform matroid, the same spectral
index, and the same Poincar\'e polynomial, but no common nonzero scalar carries
one set to the other.

Let $Z(\mathfrak s):=\{x\in\mathfrak s:[x,y]=0\text{ for all }y\in\mathfrak s\}$
denote the center of $\mathfrak s$. We first consider the center of $\mathfrak s$ and the factor $z_0$ of the
characteristic polynomial $Q_{\mathfrak s}$. As above, let $m_0$ be the
multiplicity of the zero nilradical weight, and let $\mathfrak n_0$ be the sum
of the blocks with zero diagonal weight, so $\dim\mathfrak n_0=r_0$. In light
of Theorems~\ref{thm:nw-structural-form} and~\ref{thm:nw-complex-block-form},
the strictly triangular parts of the zero-weight blocks can reduce the common
kernel even though their diagonal weights vanish. Proposition~\ref{prop:abelian-center-zero-weight}
below makes the relation between $m_0$ and the dimension of the center explicit.

Let $\iota:\mathbb{C}^r\to\mathfrak n$ be the coordinate isomorphism determined
by the basis $\{n_1,\ldots,n_r\}$. Thus $A^\alpha$ represents the linear map
$\operatorname{ad}(x_\alpha)|_{\mathfrak n}$ under $\iota$.
Set $K:=\bigcap_{\alpha=1}^k\ker A^\alpha$.

\begin{proposition}
\label{prop:abelian-center-zero-weight}
Let $\mathfrak s$ be as in Theorem~\ref{thm:nw-structural-form}. Then
\[
Z(\mathfrak s)=\iota(K)
\subseteq\mathfrak n_0,
\qquad
e_0(\mathfrak s)=k+r_0\geq k+\dim Z(\mathfrak s).
\]
The final inequality can be strict when a zero-weight block has a nonzero
strictly triangular part.
\end{proposition}

\begin{proof}
Since $Z(\mathfrak s)$ is a nilpotent ideal, the characterization of the
nilradical in a solvable complex Lie algebra gives
$Z(\mathfrak s)\subseteq\mathfrak n$. If $x\in Z(\mathfrak s)$ and
$\iota(\mathbf c)=x$, then $A^\alpha\mathbf c=0$ for $1\leq\alpha\leq k$
because $[x_\alpha,x]=0$.

Conversely, if
$\mathbf c\in\bigcap_{\alpha=1}^k\ker A^\alpha$, then
$x=\iota(\mathbf c)$ lies in $\mathfrak n$. Hence $x$ commutes with
$\mathfrak n$ because $\mathfrak n$ is abelian, and it commutes with each
$x_\alpha$ because $A^\alpha\mathbf c=0$. Therefore $x\in Z(\mathfrak s)$.
On a nonzero-weight block, at least one diagonal scalar $a_j^\alpha$ is
nonzero, so the corresponding $A^\alpha$ is invertible on that block because a
nonzero scalar plus a nilpotent matrix is invertible. Thus the common kernel
is contained in the zero-weight subspace $\mathfrak n_0$, so
$\dim Z(\mathfrak s)\leq r_0$. Finally,
Theorem~\ref{thm:abelian-block-factorization} gives
$e_0(\mathfrak s)=k+r_0$, and the stated inequality follows.
\end{proof}

For comparison with the extremal construction in
\cite{NdogmoWinternitz1994}, assume in this paragraph that $\mathfrak s$ is
nonnilpotent and indecomposable and that $1\leq k<r$. Put
$m:=r-k$, the difference between the nilradical dimension and the complement
dimension, and put $d=k(k-1)/2$. When $m>d$, write
$N_0:=r-k(k+1)/2=(k+1)q_0+\nu$, with $0\leq\nu\leq k$.
Ndogmo and Winternitz prove that the largest possible center dimension among
such algebras is
\[
 d_M=
 \begin{cases}
 m,
   &m\leq d,\\[1mm]
 \dfrac{2kr-k(k+1)}{2(k+1)},
   &m>d,\ \nu=0,\\[3mm]
 \dfrac{2kr-k(k+1)+2(\nu-k-1)}{2(k+1)},
   &m>d,\ \nu\geq1.
 \end{cases}
\]
The canonical realization begins with the $k$ weights
$\mathbf e_1,\ldots,\mathbf e_k$. Except when $\nu=1$, all remaining weights
are zero; when $\nu=1$, there is one additional nonzero weight
$\mathbf a=(a^1,\ldots,a^k)$. These cases are listed in
Table~\ref{tab:nw-maximal-center-spectra}.

In the case $\nu=1$, define
\[
\operatorname{supp}(\mathbf a):=
\{\alpha\in\{1,\ldots,k\}:a^\alpha\neq0\},
\qquad q:=|\operatorname{supp}(\mathbf a)|.
\]
Set
\[
\epsilon(\mathbf a):=
\begin{cases}
0,&\mathbf a\in\{\mathbf e_1,\ldots,\mathbf e_k\},\\
1,&\mathbf a\notin\{\mathbf e_1,\ldots,\mathbf e_k\}.
\end{cases}
\]
Finally, let $B_{\mathbf a}:=\{\mathbf e_1,\ldots,\mathbf e_k,\mathbf a\}$.

\begin{lemma}\label{lem:max-center-extra-weight}
In the $\nu=1$ case, if
$\mathbf a=\mathbf e_i$, then $\kappa=k+1$ and
$P_{\mathfrak s}(t)=(1+t)^k$. If
$\mathbf a\notin\{\mathbf e_1,\ldots,\mathbf e_k\}$ and
$\sum_\alpha a^\alpha\neq1$, then $\kappa=k+2$ and
$P_{\mathfrak s}(t)=(1+t)^{k+1}$. If
$\sum_\alpha a^\alpha=1$, then $q\geq2$, $\kappa=k+2$, and
\[
P_{\mathfrak s}(t)
=(1+t)^{k-q}
\left(\sum_{i=0}^{q-1}\binom{q+1}{i}t^i+q t^q\right).
\]
\end{lemma}

The last formula follows because the basis weights in
$\operatorname{supp}(\mathbf a)$ together with $\mathbf a$ form a circuit,
while the remaining basis weights are coloops. In particular, the formula
with $q=k$ applies when every coordinate of $\mathbf a$ is nonzero.

Table~\ref{tab:nw-maximal-center-spectra} records the resulting invariants.
The only parameter-dependent row is $\nu=1$. In that row, $\mathbf a$ may
coincide with a basis weight or satisfy an affine dependence involving the
basis weights in its support. Let $0^{[m]}$ denote the zero weight counted
with multiplicity $m$.

\begingroup
\small
\setlength{\tabcolsep}{2pt}
\Needspace{14\baselineskip}
\begin{longtable}{@{}>{\raggedright\arraybackslash}p{2.5cm}
                  >{\raggedright\arraybackslash}p{3.9cm}
                  >{\raggedright\arraybackslash}p{4.6cm}
                  >{\raggedright\arraybackslash}p{3.8cm}@{}}
\caption{Spectral invariants for maximal-center realizations.}
\label{tab:nw-maximal-center-spectra}\\
\toprule
\textbf{Case} & \textbf{Nilradical weight multiset} & \textbf{$Q_{\mathfrak s}(z_0,\mathbf{z})$} & \textbf{$\kappa$ and $P:=P_{\mathfrak s}$} \\
\midrule
\endfirsthead
\toprule
\textbf{Case} & \textbf{Nilradical weight multiset} & \textbf{$Q_{\mathfrak s}(z_0,\mathbf{z})$} & \textbf{$\kappa$ and $P_{\mathfrak s}$} \\
\midrule
\endhead
$m\leq d$
& $\{\mathbf e_1,\ldots,\mathbf e_k,0^{[m]}\}$
& $z_0^r\prod_{\alpha=1}^k(z_0+z_{r+\alpha})$
& $\kappa=k+1$, $P=(1+t)^k$ \\
$m>d$, $\nu=0$
& $\{\mathbf e_1,\ldots,\mathbf e_k,0^{[m]}\}$
& $z_0^r\prod_{\alpha=1}^k(z_0+z_{r+\alpha})$
& $\kappa=k+1$, $P=(1+t)^k$ \\
$m>d$, $\nu=1$
& $\{\mathbf e_1,\ldots,\mathbf e_k,\mathbf a,0^{[m-1]}\}$
& $z_0^{r-1}\bigl(\prod_{\alpha=1}^k(z_0+z_{r+\alpha})\bigr)
   \left(z_0+\mathbf a\mathbin{\cdot}\mathbf{z}_{\mathfrak v}\right)$
& \begin{tabular}{@{}l@{}}
$\kappa=k+1+\epsilon(\mathbf a)$,\\
$P=P_{\mathcal A_{B_{\mathbf a}}}(t)$
\end{tabular} \\
$m>d$, $2\leq\nu\leq k$
& $\{\mathbf e_1,\ldots,\mathbf e_k,0^{[m]}\}$
& $z_0^r\prod_{\alpha=1}^k(z_0+z_{r+\alpha})$
& $\kappa=k+1$, $P=(1+t)^k$ \\
\bottomrule
\end{longtable}
\endgroup

The table also makes clear that maximal center dimension is not a spectral
classification.  The Ndogmo--Winternitz classification gives alternative
maximal-center realizations
having the same value of $d_M$ but different diagonal weights.  Such
realizations have the same center dimension and can nevertheless have
different characteristic polynomials.

The two boundary cases $k=1$ and $k=r$ admit closed formulas, given in
Theorems~\ref{thm:abelian-codim-one} and
\ref{thm:abelian-minimal-nilradical}.  In the
codimension-one case, once the eigenvalues with algebraic multiplicity are
fixed, Jordan block sizes and strictly triangular entries do not affect the
characteristic polynomial.

A nonnilpotent solvable Lie algebra $\mathfrak{s}$ with $k=1$ is a one-dimensional solvable extension of its nilradical $\mathfrak{n}$. Let $A=A^1$ be the matrix defined in Theorem \ref{thm:nw-structural-form} and assume $\lambda_1,\ldots,\lambda_r$ are the eigenvalues of $A$, counting multiplicity. Set
\[
r_0:=|\{i:\lambda_i=0\}|,
\qquad
S:=\{\lambda_i:\lambda_i\neq0\}.
\]
The following result does not require the nilradical to be abelian.

\begin{theorem}
\label{thm:abelian-codim-one}
Suppose a nonnilpotent Lie algebra $\mathfrak s$ is a one-dimensional solvable extension of its nilradical. Then $S\neq\varnothing, \kappa(\mathfrak s)=1+|S|$, and $e_0(\mathfrak s)=1+r_0$. Furthermore,
\[
\begin{aligned}
Q_{\mathfrak s}(z_0,z_1,\ldots,z_{r+1})
  &=z_0\prod_{i=1}^r(z_0+\lambda_i z_{r+1}),\\
P_{\mathfrak s}(t)&=(1+t)(1+(|S|-1)t).
\end{aligned}
\]
\end{theorem}

\begin{proof}
For $k=1$, Theorem~\ref{thm:abelian-block-factorization} gives
\[
Q_{\mathfrak s}(z_0,z_1,\ldots,z_{r+1})
=z_0\det(z_0I_r+z_{r+1}A).
\]
The eigenvalues of $A$ are $\lambda_1,\ldots,\lambda_r$, counting
algebraic multiplicity, so this determinant is
$z_0\prod_i(z_0+\lambda_i z_{r+1})$.  Since $A$ is nonnilpotent,
at least one $\lambda_i$ is nonzero, so $S\neq\varnothing$.  The factor
$z_0$ from the quotient and
the factors with $\lambda_i=0$ give the stated value of $e_0$, while the
distinct nonzero eigenvalues give the stated value of $\kappa$.

The remaining arrangement consists of $|S|$ distinct lines through the
origin in $\mathbb{C}^2$.  If $|S|=1$, its complement has Poincar\'e polynomial
$1+t$, which agrees with the formula.  If $|S|\geq2$, its intersection
lattice
has one rank-zero element, $|S|$ rank-one elements, and one rank-two element.
The Orlik--Solomon formula therefore gives
\[
P_{\mathfrak s}(t)=1+|S|t+(|S|-1)t^2
=(1+t)(1+(|S|-1)t).
\]
\end{proof}

Next, we consider the case $k=r$. Let $\mathfrak a$ be the two-dimensional complex Lie algebra with basis $\{x, n\}$ satisfying the bracket relation $[x,n]=n$.

\begin{theorem}
\label{thm:abelian-minimal-nilradical}
Suppose $k=r$. Then $\mathfrak s\cong\bigoplus_{i=1}^r\mathfrak a$, and
\[
Q_{\mathfrak s}(z_0,z_1,\ldots,z_{2r})
 =z_0^r\prod_{i=1}^r(z_0+z_{r+i}),
\qquad
\kappa(\mathfrak s)=r+1,
\qquad
P_{\mathfrak s}(t)=(1+t)^r.
\]
\end{theorem}

\begin{proof}
When $k=r$, Theorem~\ref{thm:nw-complex-block-form} shows that the weights
span $\mathbb{C}^r$.  Hence at least $r$ blocks have nonzero weight.  Since the
nilradical has dimension $r$, there are exactly $r$ one-dimensional blocks.
After changing the basis of the complement, their weights are the standard
basis vectors. Identifying each matrix $A^i$ with the linear map
$\operatorname{ad}(x_i)|_{\mathfrak n}$ represented in the chosen nilradical
basis, we have $A^i(n_i)=n_i$ and $A^i(n_j)=0$ for $i\neq j$.

The Jacobi relations show that the extension vector
$R_{ij}=[x_i,x_j]$ can have components only in the directions $n_i$ and
$n_j$; when $r\geq3$, this follows by applying the relation for a triple
containing $i,j$ and a third index, and the case $r=2$ is immediate.  The
shifts $x_i\mapsto x_i+t_i$ from Theorem~\ref{thm:nw-structural-form}
remove these remaining components.  Thus the brackets between distinct
summands can be taken to vanish.  The characteristic polynomial and the
Poincar\'e polynomial then follow from the $r$ standard weights in
Theorem~\ref{thm:abelian-block-factorization}.  The resulting algebra is a
direct sum of $r$ proper ideals when $r>1$, so it is indecomposable only for
$r=1$.
\end{proof}

We next consider the weight patterns associated with the four complex
partitions for $\dim\mathfrak s=8$, $r=5$, and $k=3$.

Ndogmo and Winternitz \cite{NdogmoWinternitz1994} work out all complex partitions of the commuting action when $\dim\mathfrak s=8$, $\dim\mathfrak n=5$, and $\dim\mathfrak v=3$.  Their four partitions of $5$ are $3+1+1$, $2+2+1$, $2+1+1+1$, and $1+1+1+1+1$.
The rows record the corresponding weight multisets.  Different commuting
families, choices of nilpotent subalgebras, and extension tensors may give the
same row because they do not change the diagonal weights.  Thus Table~\ref{tab:nw-853-spectra} below 
records characteristic-polynomial data rather than compatible extension data
or isomorphism classes of Lie algebras.  The partition determines the
multiplicities of the displayed weights, and hence the exponents in
$Q_{\mathfrak s}$; the off-diagonal matrix parameters do not enter the table.
Put
$\mathbf e_1=(1,0,0)$, $\mathbf e_2=(0,1,0)$,
$\mathbf e_3=(0,0,1)$, and let
$\mathbf a=(a^1,a^2,a^3)$ and $\mathbf b=(b^1,b^2,b^3)$ denote the two
additional weights. For each row, let $\mathcal W$ denote the displayed
weight multiset and let $S$ be the set of its distinct nonzero members.
Here, generic means that the parameters avoid the exceptional values at
which weights collide or become affinely dependent.

\begingroup
\small
\setlength{\tabcolsep}{2pt}
\begin{longtable}{@{}>{\centering\arraybackslash}p{2.5cm}
                  >{\raggedright\arraybackslash}p{3.4cm}
                  >{\raggedright\arraybackslash}p{5.4cm}
                  >{\raggedright\arraybackslash}p{3.5cm}@{}}
\caption{Spectral invariants for the four partitions with $\dim\mathfrak s=8$.}
\label{tab:nw-853-spectra}\\
\toprule
\textbf{Partition} & \textbf{Weight multiset $\mathcal W$} & \textbf{$Q_{\mathfrak s}$} & \textbf{$\kappa$ and $P_{\mathfrak s}$} \\
\midrule
\endfirsthead
\toprule
\textbf{Partition} & \textbf{Weight multiset $\mathcal W$} & \textbf{$Q_{\mathfrak s}$} & \textbf{$\kappa$ and $P_{\mathfrak s}$} \\
\midrule
\endhead
\bottomrule
\endfoot
\endlastfoot

$3+1+1$
& $\{\mathbf e_1^{[3]},\mathbf e_2,\mathbf e_3\}$
& $z_0^3(z_0+z_6)^3(z_0+z_7)(z_0+z_8)$
& $\kappa=4$, $P=(1+t)^3$ \\
\addlinespace
$2+2+1$
& $\{\mathbf e_1^{[2]},\mathbf e_2^{[2]},\mathbf e_3\}$
& $z_0^3(z_0+z_6)^2(z_0+z_7)^2(z_0+z_8)$
& $\kappa=4$, $P=(1+t)^3$ \\
\addlinespace
$2+1+1+1$
& $\{\mathbf e_1^{[2]},\mathbf e_2,\mathbf e_3,\mathbf b\}$
& $z_0^3(z_0+z_6)^2(z_0+z_7)(z_0+z_8)
   (z_0+b^1z_6+b^2z_7+b^3z_8)$
& \begin{tabular}{@{}p{3.5cm}@{}}
$\kappa=1+|S|$\\
$P=P_{\mathcal A_S}(t)$\\
\emph{Generic case:}\\
$\kappa=5$, $P=(1+t)^4$
\end{tabular} \\
\addlinespace
$1+1+1+1+1$
& $\{\mathbf e_1,\mathbf e_2,\mathbf e_3,\mathbf a,\mathbf b\}$
& $z_0^3(z_0+z_6)(z_0+z_7)(z_0+z_8)
   (z_0+a^1z_6+a^2z_7+a^3z_8)
   (z_0+b^1z_6+b^2z_7+b^3z_8)$
& \begin{tabular}{@{}p{3.5cm}@{}}
$\kappa=1+|S|$\\
$P=P_{\mathcal A_S}(t)$\\
\emph{Generic case:}\\
$P=1+5t+10t^2+10t^3+4t^4$
\end{tabular} \\
\end{longtable}
\endgroup

Table~\ref{tab:nw-853-spectra} has the same interpretation.  In the third row,
$\mathbf b=0$ or $\mathbf b=\mathbf e_i$ gives
$\kappa=4$ and $P=(1+t)^3$.  If
$\mathbf b\notin\{0,\mathbf e_1,\mathbf e_2,\mathbf e_3\}$, then
$b^1+b^2+b^3\neq1$ gives the \emph{generic} value
$\kappa=5$ and $P=(1+t)^4$.  If
$b^1+b^2+b^3=1$ and all three coordinates of $\mathbf b$ are nonzero, then
$\kappa=5$ and $P=1+4t+6t^2+3t^3$.  If exactly one coordinate of
$\mathbf b$ is zero, then a three-element circuit occurs and
$P=1+4t+5t^2+2t^3$.

The parameter-dependent rows remain valid when weights coincide or become zero.
To compute the spectral matroid and Poincar\'e polynomial, delete zero
weights and retain only one copy of each distinct nonzero weight before forming
the coefficient vectors $(1,\mathbf b)$.  For the characteristic polynomial, the exponent
$e_0$, and spectral equivalence, all weight multiplicities, including the
zero-weight multiplicity, must be retained.  In the fourth row, five distinct
nonzero weights have \emph{generic} value
$P=1+5t+10t^2+10t^3+4t^4$.  For this reduced set, the spectral matroid is
obtained by computing the ranks of all row-submatrices
of the coefficient matrix, whose rows are
$(1,\mathbf e_1)$, $(1,\mathbf e_2)$, $(1,\mathbf e_3)$,
$(1,\mathbf a)$, and $(1,\mathbf b)$.  If every three-element subset is
already known to be affinely independent, then the $4\times4$ minors
determine which four-element subsets are dependent.  At special parameter
values, three-element circuits must also be detected using the corresponding
lower-order minors.  For example,
$\mathbf a=(\tfrac12,\tfrac12,0)$ gives, with
$\widehat{\mathbf u}:=(1,\mathbf u)$, the affine relation
$\widehat{\mathbf a}=\tfrac12\widehat{\mathbf e}_1
+\tfrac12\widehat{\mathbf e}_2$.
Thus the table entry $P_{\mathcal A_S}(t)$ applies at every parameter value.

\subsection{A Dimension 9 Worked Example}
\label{sec:abelian-worked-example}

We further illustrate the preceding results with a dimension-$9$ example. We compute its characteristic polynomial, identify the spectral matroid, and compare it with a second algebra
having the same characteristic polynomial but a different derived algebra.

Consider $\mathfrak{s}:=\mathfrak{n}\oplus \mathfrak{v}$, where $\mathfrak n=\operatorname{span}\{n_1,\ldots,n_6\}$ is abelian, and 
$\mathfrak v=\operatorname{span}\{x_1,x_2,x_3\}$.
We use the ordered basis
$\mathcal B=(n_1,n_2,n_3,n_4,n_5,n_6,x_1,x_2,x_3)$ and define the action of
$\mathfrak v$ on $\mathfrak n$ by the diagonal matrices
\[
A^1=\operatorname{diag}(1,0,0,2,7,0),\quad
A^2=\operatorname{diag}(0,1,0,3,11,0),\quad
A^3=\operatorname{diag}(0,0,1,5,13,0).
\]
More explicitly, set $[x_\alpha,u]=A^\alpha u$ for $u\in\mathfrak n$ and
set
\[
[x_1,x_2]=n_6,
\qquad [x_1,x_3]=[x_2,x_3]=0.
\]
All other brackets are determined by skew-symmetry or are zero. Since $\mathfrak n$ is abelian, the only
Jacobi identity involving the extension bracket that is not immediate is the
one for $(x_1,x_2,x_3)$. Its three terms vanish: the term involving the
extension bracket is
\[
[x_3,[x_1,x_2]]=[x_3,n_6]=0.
\]
The other two terms contain $[x_2,x_3]$ or $[x_3,x_1]$, so they are zero by
definition.
Thus these brackets define a Lie algebra.

The diagonal weights are
\[
\begin{aligned}
\alpha_1&=(1,0,0),&\alpha_2&=(0,1,0),&\alpha_3&=(0,0,1),\\
\alpha_4&=(2,3,5),&\alpha_5&=(7,11,13),&\alpha_6&=0.
\end{aligned}
\]

These weights also help identify the center and nilradical. If
$u+c_1x_1+c_2x_2+c_3x_3$ is central, commuting with
$n_1,n_2,n_3$ gives $c_1=c_2=c_3=0$. The common kernel of
$A^1,A^2,A^3$ is then $\mathbb{C}n_6$, so
$Z(\mathfrak s)=\mathbb{C}n_6$. On the other hand, the first three diagonal
entries of $c_1A^1+c_2A^2+c_3A^3$ are $c_1,c_2,c_3$. A nilpotent combination
must therefore be zero, so the matrices are linearly nilindependent. Hence
$\mathfrak n$ is the nilradical.

For
$y=z_1n_1+\cdots+z_6n_6+z_7x_1+z_8x_2+z_9x_3$, the adjoint matrix has the
block form
\[
\operatorname{ad}(y)=
\begin{pmatrix}
D(z_7,z_8,z_9)&*\\
0&0_3
\end{pmatrix},
\qquad
D=\operatorname{diag}(z_7,z_8,z_9,
2z_7+3z_8+5z_9,7z_7+11z_8+13z_9,0).
\]
The entries in the upper-right block contain the nilradical coordinates and the
extension bracket, but they do not affect the determinant. The lower-right
block is zero because $\mathfrak s/\mathfrak n$ is abelian. Therefore
\[
\begin{aligned}
Q_{\mathfrak s}(z_0,z_1,\ldots,z_9)
&=z_0^3\det(z_0I_6+D)\\
&=z_0^4(z_0+z_7)(z_0+z_8)(z_0+z_9)\\
&\qquad\cdot(z_0+2z_7+3z_8+5z_9)(z_0+7z_7+11z_8+13z_9).
\end{aligned}
\]
The quotient contributes $z_0^3$ and the zero weight contributes one more
copy, so $e_0=4$ and $\kappa=6$.

Remove the factor $z_0^4$. The five remaining factors are
\[
\begin{aligned}
L_1&=z_0+z_7, &
L_2&=z_0+z_8, &
L_3&=z_0+z_9,\\
L_4&=z_0+2z_7+3z_8+5z_9, &
L_5&=z_0+7z_7+11z_8+13z_9.
\end{aligned}
\]
Their coefficient matrix is
\[
M=
\begin{pmatrix}
1&1&0&0\\
1&0&1&0\\
1&0&0&1\\
1&2&3&5\\
1&7&11&13
\end{pmatrix}.
\]
For $1\leq i\leq5$, let $M_{\widehat{i}}$ be the $4\times4$ matrix obtained
by deleting the $i$th row of $M$. For example,
\[
M_{\widehat{1}}=
\begin{pmatrix}
1&0&1&0\\
1&0&0&1\\
1&2&3&5\\
1&7&11&13
\end{pmatrix},
\qquad
\det M_{\widehat{1}}=-3.
\]
The five determinants are
\[
\bigl(\det M_{\widehat{1}},\ldots,\det M_{\widehat{5}}\bigr)
=(-3,\ 9,\ 33,\ 30,\ 9).
\]
All five are nonzero, so every four of the five coefficient vectors is
independent. Since there are only four coordinates, the spectral matroid is
therefore the uniform matroid $U_{4,5}$. Its unique circuit gives the relation
\[
-L_1-3L_2+11L_3-10L_4+3L_5=0.
\]

For $U_{4,5}$, the Whitney numbers are
$(w_0,w_1,w_2,w_3,w_4)=(1,5,10,10,4)$. Hence
\[
P_{\mathfrak s}(t)=1+5t+10t^2+10t^3+4t^4,
\qquad (b_0,b_1,b_2,b_3,b_4)=(1,5,10,10,4).
\]

Let $\mathfrak s_0$ be the solvable extension of $\mathfrak{n}$ obtained by replacing $[x_1,x_2]=n_6$ with $[x_1,x_2]=0$, while keeping the matrices $A^1,A^2,A^3$ unchanged. The
extension bracket appears only in the upper-right block of $\operatorname{ad}(y)$,
so the characteristic polynomial is unchanged:
\[
Q_{\mathfrak s_0}(z_0,z_1,\ldots,z_9)
=Q_{\mathfrak s}(z_0,z_1,\ldots,z_9).
\]
Their derived algebras, that is, their commutator subalgebras, are different:
\[
[\mathfrak s_0,\mathfrak s_0]
=\operatorname{span}\{n_1,\ldots,n_5\},
\qquad
[\mathfrak s,\mathfrak s]
=\operatorname{span}\{n_1,\ldots,n_6\}.
\]
Their derived algebras have dimensions $5$ and $6$, respectively, so
$\mathfrak s_0\not\cong\mathfrak s$, although their characteristic
polynomials are identical.

Finally, change the weight of $n_5$ from $(7,11,13)$ to $(2,3,5)$. It then
coincides with the weight of $n_4$, so $L_5=L_4$ and the last two factors merge:
\[
Q_{\mathfrak s}(z_0,z_1,\ldots,z_9)
=z_0^4(z_0+z_7)(z_0+z_8)(z_0+z_9)
  (z_0+2z_7+3z_8+5z_9)^2.
\]
The number of distinct factors other than $z_0$ drops from five to four, so
the spectral index $\kappa$ drops from $6$ to $5$; the reduced arrangement now has four
independent factors and $P_{\mathfrak s}(t)=(1+t)^4$. Thus a repeated weight
changes the spectral matroid and the topology even though the
calculation still uses the same block form.

\section{Concluding Remarks}

The classification of finite-dimensional solvable Lie algebras over
$\mathbb{C}$ is a longstanding open problem, with isomorphism classes
proliferating rapidly even in low dimensions. In this paper, we developed a
spectral approach based on the multivariate invariant
$Q_{\mathfrak{s}}(z_0,\mathbf{z})=\det\bigl(z_0I+\sum_i z_i\operatorname{ad} x_i\bigr)$.

The main application is the spectral analysis of solvable complex Lie algebras
with abelian nilradical, including their classification up to spectral
equivalence through the corresponding weight multisets.  The
Ndogmo--Winternitz normal form makes the general factorization and arrangement
theory explicit: triangular block sizes give
the weight multiplicities, Proposition~\ref{prop:abelian-weight-realization}
realizes every spanning collection of nonzero weights, and
Remark~\ref{cor:nw-spectral-equivalence} and
Corollary~\ref{cor:abelian-affine-matroid} identify the resulting spectral equivalence
criterion and matroid.
The center bound and the tables apply these statements to the families
treated by Ndogmo and Winternitz, while the worked example exhibits two
non-isomorphic extension choices with the same characteristic polynomial.
This makes the abelian-nilradical class a useful test case for the distinction
between isomorphism and spectral equivalence.

The subsequent papers apply these results to several classes of nilradicals.
In each case, the defining brackets impose additional relations among the
weights or determine them from the classified extension derivations.  The
resulting spectral matroids and parameter spaces are different, but the same
factorization makes them comparable.

Several directions remain open.

\begin{enumerate}
\item \emph{Other nilradicals.}
Which nilradicals impose intrinsic restrictions on the associated spectral
matroids? Such restrictions could give meaningful classification results even
when a complete classification up to isomorphism is out of reach.

\item \emph{Spectral rigidity.}
Within a parameterized family, when does the nilradical weight multiset
determine the Lie algebra up to isomorphism? Answering this would identify
when spectral equivalence agrees with isomorphism and which extension
data are not determined by the characteristic polynomial.

\item \emph{Geometric applications.}
Given a finite-dimensional real Lie algebra $\mathfrak g$, let $G$ be the
corresponding simply connected real Lie group and let $\mathfrak g_{\mathbb C}$
denote its complexification. It is natural to ask what the spectral data of
$\mathfrak g_{\mathbb C}$ say about geometric features of $G$ and its
homogeneous spaces.
\end{enumerate}

These results place spectral equivalence between isomorphism and comparison by
coarser numerical invariants.  In the abelian-nilradical case, a commuting
matrix pencil gives the weights and their multiplicities directly, while the
linear-dependence relations among the coefficient vectors $(1,\alpha)$ determine
the associated arrangement matroid.  The Ndogmo--Winternitz normal form also makes explicit
which extension parameters are omitted by the characteristic polynomial.
Spectral equivalence therefore provides a computable comparison for
parameterized families which is coarser than isomorphism.

\section*{Acknowledgements}

This project was supported by Innovation Funding from the Minerva Center for High-Impact Practices and SUNY System Administration. The second author is supported partially by the Simons Foundation Grant No. MP-TSM-00002315.

\end{document}